\documentclass[12pt]{article}

\usepackage{epic}
\usepackage{amsmath}[1996/11/01]
\usepackage{url,amssymb,amsthm,amsfonts,latexsym,epsfig,graphics}
 \def\beql#1#2\eeql{\begin{equation}\label{#1}#2\end{equation}}

\newcommand{\knubbel}{\begin{picture}(5,5)(3,3)  \put(1,6){\circle*{4}} \end{picture}}

\DeclareMathOperator{\ind}{ind}

\DeclareMathOperator{\Fix}{Fix}

\DeclareMathOperator{\Br}{Br}
\DeclareMathOperator{\Gal}{Gal}
\DeclareMathOperator{\Aut}{Aut}

\DeclareMathOperator{\Kappa}{{\mathfrak K}}
\DeclareMathOperator{\Rho}{{\mathfrak R}}
\DeclareMathOperator{\Chi}{{\mathfrak X}}

\DeclareMathOperator{\Irr}{Irr}
\DeclareMathOperator{\End}{End}

\DeclareMathOperator{\disc}{disc}
\DeclareMathOperator{\GL}{GL}

\DeclareMathOperator{\diag}{diag}

\DeclareMathOperator{\SU}{SU}
\DeclareMathOperator{\id}{id}

\DeclareMathOperator{\SL}{SL}

\theoremstyle{plain}
\newtheorem{theorem}{Theorem}
\newtheorem{lemma}[theorem]{Lemma}

\newtheorem{proposition}[theorem]{Proposition}
\newtheorem{corollary}[theorem]{Corollary}
\newtheorem{definition}[theorem]{Definition}
\theoremstyle{remark}

\newtheorem{remark}[theorem]{Remark}
\newtheorem{example}[theorem]{Example}

\numberwithin{theorem}{section}

\newcommand{\Z}{{\mathbb{Z}}}
\newcommand{\Q}{{\mathbb{Q}}}
\newcommand{\F}{{\mathbb{F}}}

\newcommand{\N}{{\mathbb{N}}}
\newcommand{\R}{{\mathbb{R}}}

\newcommand{\C}{{\mathbb{C}}}
\newcommand{\trace}{\mbox{trace}}

\begin{document}
\title{Unitary discriminants of characters}
 \author{Gabriele Nebe\footnote{
 Lehrstuhl f\"ur Algebra und Zahlentheorie, 
 RWTH Aachen University, Germany, gabriele.nebe@rwth-aachen.de}}

\maketitle

 \begin{center}
	 {\sl In memory of Richard Parker}
 \end{center}

 {\sc Abstract.} Together with David Schlang we computed the 
 discriminants of the invariant Hermitian forms for all indicator $o$ 
 even degree irreducible characters of ATLAS groups 
 supplementing the tables of orthogonal determinants 
 computed in collaboration with Richard 
 Parker, Tobias Braun and Thomas Breuer. 
 The methods that are used in the unitary case are described in this
 paper. 
 An ordinary character has a well defined unitary discriminant algebra,
 if and only if
 it is unitary stable, i.e. all irreducible 
 unitary constituents have even degree. 
Computations for large degree characters are only possible because  
 of a new method called {\em unitary condensation}. 
Often, a suitable automorphism helps to single out a square class of the 
real subfield of the character field consisting of representatives of 
the discriminant of the invariant Hermitian forms. 
This square class can then be determined 
modulo enough primes.
     \\
  {\sc keywords:}  unitary representations of finite groups;  
unitary discriminants; invariants of Hermitian forms; Clifford invariants of quadratic forms; discriminant algebras of involutions of the second kind

\section{Introduction} 

Let $K$ be a field, $G$ a finite group and $n\in \N$.
A group homomorphism $\rho : G\to \GL_n(K)$  
is called a $K$-representation of $G$. 
 Often $\rho (G)$ is contained in a smaller classical group, 
 such as symplectic, unitary or orthogonal groups. 
 This paper 
 continues a long term project of the author with 
 Richard Parker to specify these classical groups
 for finite fields and number fields. 
 The orthogonal groups are dealt with in the joint paper 
 \cite{OrthogonalStability} and symplectic groups are 
 essentially uniquely determined by the irreducible character.
 This paper 
 approaches computational challenges to precise 
 unitary representations over number fields.
Whereas over finite fields and also over the real numbers 
there is a unique compact unitary group of any given dimension, 
there are infinitely many such groups over number fields if the dimension 
of the underlying Hermitian space is even.
These unitary groups are uniquely determined by 
the discriminant of the invariant definite Hermitian forms.

To be more precise 
let $\chi $ be an ordinary indicator $o$ irreducible
 character of the finite group $G$ 
 and put $L := \Q(\chi )$ to denote its 
character field, i.e. the number field that is generated by 
the character values of $\chi $. Then $L$ is a totally complex number field. 
Let $K$ denote its real subfield and $N:=N_{L/K}(L^{\times}) $ the 
norm subgroup of $K^{\times }$.  Any representation $\rho $ affording the 
 character $\chi $ fixes definite Hermitian forms. 
All these invariant forms have the same discriminant if and only if 
 the degree $\chi (1)$ is even. Then we say that $\chi$ is 
 {\em unitary stable}.
 
 Let $\chi $ be such an ordinary irreducible even degree indicator
 $o$ character. 
Due to Schur indices there might not be a representation $\rho $ over 
the character field $L$ affording $\chi$. 
Nevertheless the character defines a central simple 
$L$-algebra with an involution of the second kind, for which 
  \cite[Section 10]{BookOfInvolutions} defines a 
 Brauer class $[D]$ 
 of a quaternion algebra over $K$, the so called {\em discriminant algebra}, 
 that allows to read off the discriminant of all $\rho (G) $-invariant
 Hermitian forms for any representation $\rho $ of $G$ affording 
 the character $\chi $ (cf. Remark \ref{discalggen}). 
 If  $\chi $ is quasi-split (Definition \ref{defdisc}) 
 (i.e. the quotient of degree and Schur index is even)
 then $L$ splits $D$ and hence 
 $D= (L,\delta )_K$ for a well defined $\delta \in K^{\times}/N$, 
 which we then call the {\em unitary discriminant} of $\chi $
 (see Definitions \ref{defdisc} and \ref{defdiscrchi}). 
In our heavy computations we could always avoid to work with discriminant algebras; 
the relevant indicator $o$ characters all 
turned out to have Schur index $1$, i.e. they are characters of a representation over 
the character field. 
However, for the sake of future applications, we give the basics for 
discriminant algebras and formulate most results in these terms. 
Thanks to Remark \ref{discalggen} the reader may replace the discriminant algebra 
(which is an element of the Brauer group of $K$) 
by the discriminant  (represented by a non-zero element of $K$).

The computational methods,
developed for orthogonal characters in \cite{OrthogonalStability},
 such as restriction,
 induction, tensor products and symmetrizations can also be used 
 in the unitary case.
 Also the decomposition number techniques 
can be varied to obtain restrictions for the unitary discriminants
(see Section \ref{decomposition}),
for instance showing that all primes ramifying in the 
discriminant algebra also do divide the group order 
(see Corollary \ref{ramdiv}). 

  However there is a great difference to the orthogonal case: 
  The discriminant of a given $L/K$-Hermitian form $H$ is only well defined 
  modulo norms, $\disc(H) \in K^{\times }/N_{L/K}(L^{\times })  $. 
For finite fields, the norm map is surjective, so
it seems to be impossible to obtain any information about 
$\disc (H)$ by reducing the representation modulo primes. 
The latter technique is important in \cite{OrthogonalStability} 
to handle large degree characters, 
for which we computed composition factors of a condensed module 
(see Section \ref{orthcond}) 
modulo some well chosen primes (usually not dividing the group order) 
to get enough information on the square class containing the
orthogonal discriminant. 
It is amazing that also this condensation method can be transferred 
to the unitary case: 
If there exists a suitable automorphism $\alpha $ of order 2 of $G$, 
interchanging the indicator $o$ character $\chi $ with 
its complex conjugate $\overline{\chi} = \chi \circ \alpha $, 
then there is a square class of the real subfield of the 
character field of $\chi $, called {\em the $\alpha $-discriminant of $\chi $},
that relates 
to the unitary discriminant of $\chi $ (see Theorem \ref{fixalpha}). 
In complete analogy to the orthogonal case, we then can 
use skew adjoint units in the $\alpha $-fixed algebra 
to compute the $\alpha $-discriminant of $\chi $ (see Section \ref{unitarycond}). 

We start by recalling the basic notions for quadratic and Hermitian forms. 
As we are mainly interested in totally positive definite 
Hermitian forms over quadratic 
complex extensions $L$ of totally real number fields  $K$,
Section \ref{NUMBER} gives the important facts for this situation. 
In particular Hasse's norm theorem gives criteria for 
computing the norm subgroup $N_{L/K}(L^{\times })$ of $K^{\times }$ 
from local data. 
Invariant lattices provide a tool to use modular reduction 
to restrict the possible unitary discriminants. 
Therefore some facts about lattices in local Hermitian spaces are
recalled in Section \ref{locallat}.

Section \ref{alginvo} recalls and extends the relevant 
notions for discriminants of 
algebras with involution. Many of the results can be found 
in \cite[Section 10]{BookOfInvolutions}.

The next section introduces unitary discriminants and discriminant algebras of characters.
To appropriately use modular reduction techniques, 
we need to have a finer notion of indicator for 
 irreducible Brauer characters: 
There are two different sorts of indicator $o$ Brauer characters;
those, for which the corresponding simple module carries a non-zero 
invariant Hermitian form, and those for which there is no such invariant
form. The former characters are called {\em unitary}, as the corresponding
representation embeds into the unitary group. 
To simplify notation we also say that  all real Brauer characters and all 
ordinary characters are unitary.
Then a character is called {\em unitary stable} 
(see Definition \ref{unitarystable}) if and only if all its 
unitary constituents have even degree. 
Unitary stable quasi-split ordinary characters 
are those that have a well defined unitary
discriminant.
Section \ref{methods} collects important character theoretic methods that
can be automatised to compute and validate unitary discriminants. 
  For orthogonal characters there was no need to deal with covering
  groups, as there the orthogonal discriminants are predicted by 
  theoretical methods. 
 For unitary characters the 
  $Q_8$-trick described in Section \ref{2G} is a shortcut to obtain
  unitary discriminants of 
 faithful characters of certain even covering groups.

The next section lays the ground for applying condensation techniques 
to compute unitary discriminants. 
The discriminant of a symmetric bilinear form can be computed as 
the discriminant of its adjoint involution, and hence as the 
square class of the discriminant of any skew-adjoint unit. 
The adjoint involution of an Hermitian form $H$ is an involution of
second kind. Also here we can obtain the discriminant of $H$
intrinsically from the algebra $\End(H)$
with involution (see Remark \ref{discalggen}).
There are two obstacles, however:
The discriminant algebra from Remark \ref{discalggen} is hard to determine
in general and the discriminant of $H$ is only defined modulo norms.
Orthogonal subalgebras (Definition \ref{orthsub}), if they exist, 
provide a workaround.
Orthogonal subalgebras can be obtained as $\alpha $-fixed algebras 
(Definition \ref{deforthsub}) for certain outer automorphisms $\alpha $. 
Section \ref{condensation} shows how this can be applied 
to compute unitary discriminants using modular condensation techniques. 
After reporting on highly non-trivial computations for the Harada-Norton group,
 we give a small example for $U_3(7)$ which on the one hand 
illustrates the condensation techniques and on the other hand 
generalises to compute the 
unitary discriminants of the groups $U_3(q)$ for odd prime powers $q$
(see \cite{SL3SU3unitary})
(and very likely for other finite groups of Lie type).
The last section is devoted to report on one interesting example, 
the group $3.ON$, the 
Schur cover of the sporadic simple O'Nan group. 
Besides giving examples for the various methods described in this paper,
the results are interesting, as this is one of the rare cases, where 
some unitary discriminants are not rational. 
Non-rational orthogonal discriminants occur much more frequently 
(see \cite{OD}).

Special thanks go to David Schlang, who compiled a large list of 
unitary discriminants during his master thesis \cite{MasterSchlang} 
and his fellowship awarded by the SFB TRR 195, Computeralgebra.
With his help we obtained the unitary discriminants of 
 all indicator $o$ even degree ordinary characters 
 for all groups in the 
 ATLAS of finite groups \cite{ATLAS} 
  of order smaller or equal to the order of the Harada-Norton 
  group, the current status of the database of 
  orthogonal discriminants.
I also thank Thomas Breuer, who is permanently updating the 
OSCAR database of orthogonal and unitary discriminants of characters 
\cite{OD}.

This paper is a contribution to
Project-ID 286237555 – TRR 195 – by the Deutsche Forschungsgemeinschaft
(DFG, German Research Foundation).

\section{Quadratic and Hermitian forms} 

Let $L$ be a field of characteristic $\neq 2$ and 
$\sigma $ an automorphism of $L$ of order $1$ or $2$. 
 Put $K:= \Fix _{\sigma }(L) $ to be the fixed field of 
 $\sigma $ in $L$. 
 Then either $\sigma = \id $ and $L=K$ and we talk 
 about quadratic forms or $\sigma $ has order 2
 and $[L:K] =2$. 
A Hermitian space, or, more precisely, an $L/K$-Hermitian space 
is a finite dimensional vector space $V$ over $L$ with a non-degenerate
Hermitian form $H:V\times V \to L$, i.e. a $K$-bilinear map
such that $H(av,w) = aH(v,w)$ and $H(v,w)= \sigma(H(w,v))$ for
all $v,w\in V$ and $a\in L$. 
Let 
 $$N:= N_{L/K}(L^{\times }) = \{ a \sigma(a) \mid a\in L^{\times} \}$$ 
 denote the norm subgroup of $K^{\times }$.
 Then $(K^{\times })^2 \leq N \leq K^{\times } $ and $N=(K^{\times })^2$ 
 if $L= K$.

 \begin{definition}\label{disc}
	The {\em discriminant} of $H$ is the signed determinant 
	 $$ \disc(H) := (-1)^{{{n}\choose{2}}}  \det (H_B) N \in K^{\times } / N$$ 
	where $n:=\dim_L(V)$ is the dimension of $V$ and 
	$H_B := (H(b_i,b_j))_{i,j =1 }^n \in L^{n\times n}$ is the Gram matrix of 
	$H$ with respect to any $L$-basis $B=(b_1,\ldots , b_n)$ of $V$.
\end{definition}

Of course the discriminant is already determined by the dimension and the 
determinant, so discriminant and determinant are equivalent invariants. 
However, the discriminant seems to be  the more appropriate choice, 
mostly because discriminants of hyperbolic spaces are 1. 

 For quadratic forms $H$, i.e. $L=K$, the square class $\disc(H)=d(K^{\times })^2$ 
 defines a unique \'etale $K$-algebra $\mathcal{D}(H) = K[X]/(X^2-d)$ 
 of degree 2 over $K$, the {\em discriminant algebra} of the 
 quadratic form (see \cite[Definition 3.1, Remark 3.2]{OrthogonalStability} 
 also for the correct 
 definition for fields of characteristic $2$).
For $L\neq K$ the class of the discriminant 
defines a unique quaternion algebra over $K$ which 
we call the discriminant algebra of $H$: 

 \begin{definition} 
For $a,b\in K^{\times }$ let $(a,b)_{K} $ denote the central simple
        $K$-algebra with $K$-basis $(1,i,j,ij)$ such that
        $i^2=a$, $j^2=b$, $ij=-ji$.
        If $L= K[\sqrt{\delta }]$ is a quadratic extension of $K$ we also put
        $$(L,b)_K  := (\delta,b)_{K} .$$
      For $\sigma \neq \id$   the {\em discriminant algebra} of the Hermitian
	 form $H$ with 
	 $\disc(H) = dN_{L/K}(L^{\times }) $
        is defined as the class 
	 $$\Delta (H) := [(L,d)_{K}] \in \Br_2(L,K) $$
	 where 
	 $\Br_2(L,K)$ is 
	  the exponent $2$-subgroup of the 
	 Brauer group $\Br(K)$ of central simple $K$-algebras 
	  that are split by $L$.
\end{definition}

For $a,b\in K^{\times }$ it holds that 
$(L,a)_K \cong (L,b)_K$ if and only if $aN_{L/K}(L^{\times }) 
=bN_{L/K}(L^{\times }) $.

\begin{definition} \label{Ldisc}
	For $[D] \in \Br_2(L,K)$ the {\em $L$-discriminant} of $[D]$ is defined 
	as 
	$$ \disc _L([D])  = b N_{L/K}(L^{\times }) \in K^{\times }/N_{L/K}(L^{\times })$$ 
	where $b\in K^{\times }$ is such that 
	$[D] = [(L,b)_K] $.
\end{definition}

For an $L/K$-Hermitian form $H$ the $L$-discriminant of $\Delta (H)$ 
is equal to the discriminant of $H$.

\section{Hermitian forms over number fields} \label{NUMBER}

We now assume that $L$ is a complex number field with 
totally real subfield $K:=L^+$, in particular $[L:K] = 2$ and 
$\id \neq \sigma \in \Aut_K(L)$ is the non-trivial 
$K$-linear field automorphism of $L$. 
Let $(V,H)$ be a totally positive definite Hermitian space of dimension, say, 
 $n:=\dim_L(V)$. 

Restriction of scalars turns $V$ into a 
vector space $V_K$ of dimension $2n$ over $K$. 
The Hermitian form $H$ defines a quadratic form $Q_H$ on $V_K$ 
by 
$$Q_H(v) := H(v,v)  \mbox{ for all } v\in V_K .$$ 
The following proposition is given in \cite{Scharlau} for general global fields.
Note that in \cite[Remark (10.1.4)]{Scharlau} the determinant of the
Hermitian form has to be replaced by its discriminant to obtain a correct 
statement. 

\begin{proposition}(see \cite[Remark (10.1.4), Theorem (10.1.7)]{Scharlau})
	Two Hermitian spaces $(V,H)$ and $(V',H')$ are isometric if and 
	only if the quadratic spaces $(V_K,Q_H)$ and $(V'_K,Q_{H'})$ 
	are isometric. 
	If $L=K[\sqrt{\delta}] $ then $\disc(Q_H) = \delta ^n$ and 
the Clifford invariant of $Q_H$ is  $\Delta(H)$.
\end{proposition}

As Clifford invariant, discriminant, dimension, and signatures at all 
real places of $K$ determine the isometry class of a quadratic form 
over $K$, we get the following corollary. 

\begin{corollary}
Two totally positive definite $L/K$ Hermitian forms are isometric,
if and only if they have the same dimension and the same discriminant.  
\end{corollary}

An important feature for number fields $K$ is the 
so called Hasse principle or local-global principle.
The completions of $K$ are given by the embeddings of $K$ into
the real or complex numbers (the infinite places) 
and the finite places, which are in bijection  
n to the maximal ideals of its ring of integers $\Z_K $.

\begin{itemize} 
\item 
Two $L/K$ Hermitian forms are isometric over $K$ if and only if 
they are isometric over all its completions. 
\item 
	Two quaternion algebras over $K$ are isomorphic, 
	if and only if they ramify at the same places of $K$. 
	Moreover the number of ramified places is always even. 
	Recall that  a place of $K$ is ramified in the quaternion 
	algebra $Q$, if the corresponding completion 
	of $Q$ is a division algebra. 
\item Hasse Norm Theorem for quadratic extensions:
	For $a\in K$ it holds that 
 $a\in N_{L/K} (L) $ if and only
	if $a\in N_{L_{\wp}/K_{\wp} } (L_{\wp })$ for all 
	places $\wp $ of $K$.
\end{itemize}

Let $\wp \unlhd \Z_K$ be a maximal ideal. 
Then there are three cases: 
\begin{itemize}
\item[a)]  $\wp \Z_L $ is a maximal ideal of $\Z_L$, i.e. $\wp $ is
	inert in $L/K$ and $L_{\wp} $ is the unique unramified quadratic extension of $K_{\wp }$.
\item[b)]
 $\wp \Z_L = P_1 P_2 $ is a product of two distinct 
		prime ideals $P_1=\sigma(P_2)$ of $\Z_L$, i.e. 
		$\wp $ is split in $L/K$ and $L_{\wp } \cong K_{\wp } \oplus 
		K_{\wp }$.
\item[c)]
$\wp \Z_L = P^2$ is the square of a prime ideal of $\Z_L$. 
Then $\wp $ is ramified in $L/K$.
\end{itemize}

The following
local criterion for being a norm is important for the computations.

\begin{lemma} \label{norm}
	An element $0\neq a \in K$ is a norm, $a\in N_{L/K}(L^{\times })$ if 
	and only if the following three criteria are satisfied:
	\begin{itemize}
		\item[(i)] $a$ is totally positive 
		\item[(ii)] $\nu _{\wp}(a) $ is even for all prime ideals 
			$\wp $ of $K$ that are inert in $L/K$ 
		\item[(iii)] $a \in N_{L_{\wp}/K_{\wp} } (L_{\wp }) $ 
			for all ramified primes $\wp $ of $K$.
	\end{itemize}
\end{lemma} 

\begin{proof}
	The first criterion guarantees that $a$ is a norm 
	for all infinite places of $K$, the second one that 
	$a$ is a local norm for all inert places and the 
	third one for the ramified places. For the split places 
	$L_{\wp } = K_{\wp } \oplus K_{\wp }$ and all elements of
	$K_{\wp }$ are norms. 
\end{proof}

For the real completions the signature 
parameterizes Hermitian forms. As  $H$ is totally
positive definite, this specifies $H$ for the real completions of $K$.

\begin{remark} \label{infty}
The infinite places of $K$ are ramified in the discriminant algebra 
	$\Delta (H)$ 
if and only if $\disc(H) $ consists of totally negative elements, 
i.e. $\dim _L(V) \equiv $ $2$ or $3$ mod $4$.
\end{remark} 

Let $\wp $ be a prime ideal of $K$ that is split in $L/K$.
Then the completion 
$\Delta (H)_{\wp }$ contains the completion $L_{\wp } \cong K_{\wp }
\oplus K_{\wp }$ which contains zero divisors. So
$\wp $ is not ramified in $\Delta (H)$.

\begin{remark} \label{split}
No split prime ideal ramifies in the discriminant algebra of $H$.
\end{remark}

\begin{proposition}
Let $\wp $ be a prime ideal of $K$ that is inert in $L/K$.
Then $\wp $ is ramified in $\Delta (H) $ if and only if the
$\wp $-adic valuation $\nu _{\wp } (d)$ is odd for all $d\in \disc(H)$.
\end{proposition}

\begin{proof}
Then $L_{\wp }/K_{\wp }$ is the unique unramified extension of
degree 2 of $K_{\wp }$. In particular the norm group
$$N_{L_{\wp}/K_{\wp}} = \{ x\in K_\wp \mid \nu _{\wp } (x) \in 2\Z \} .$$
	Now $\wp $ is ramified in $\Delta (H) = [(L,d)_K]$ if and only if
$d$ is not a norm in $L_{\wp }/K_{\wp }$, i.e. if and only if 
	$\nu _{\wp}(d) $ is odd. 
\end{proof}

\begin{remark} 
For a ramified place $\wp = P^2$ we can always find a representative
$d \in \Z_K$ of $\disc(H)$ such that $d\not\in \wp $ 
	(see Proposition \ref{duallat} (c)). 
Assume that $\wp $ is not a dyadic prime, i.e. $2\not\in \wp $. 
Then $\wp $ ramifies in the discriminant algebra of $H$ if and only 
if $d$ is not a square modulo $\wp $.
\end{remark} 

\subsection{Lattices}

Assume that $L/K$ is a quadratic extension of local or global number fields
and $(V,H)$ a finite dimensional non-degenerate $L/K$-Hermitian space.

\begin{definition}
	A $\Z_L$-lattice $\Lambda $ in $V$ is a finitely generated $\Z_L$-submodule of $V$ that contains a basis of $V$. 
	The {\em dual lattice} is
	$$\Lambda ^* :=\{ v\in V \mid H(v,\Lambda ) \subseteq \Z_L \} $$
	and again a $\Z_L $-lattice in $V$.
	We say that $\Lambda $ is {\em unimodular} if $\Lambda =\Lambda ^*$.
\end{definition}

\begin{example}
	Let $K=\Q $ and $L=\Q[\sqrt{-10}]$.
	Then $P_5:=(5,\sqrt{-10})  \unlhd \Z_L $ is 
	not a principal ideal. 
	Let $(V,H)$ be an $n$-dimensional $L/K$-Hermitian 
	space with orthonormal basis $B=(b_1,\ldots , b_n)$, so $H_B=I_n$ 
	and $(V,H') $ such a space with $H'_{B} = \diag(I_{n-1},1/5)$.
	Then 
	$$\Lambda = \langle B \rangle _{\Z_L}  \mbox{ and } 
	\Lambda ' := \bigoplus_{i=1}^{n-1} \Z_L b_i \oplus P_5 b_n $$ 
	are unimodular lattices in $(V,H)$ and  $(V,H')$ respectively. 
	If $n$ is even, then 
	$$ \Delta(H) = \left\{ \begin{array}{lll}  {[}(-10,-1)_{\Q }{]} = {[}(-1,-1)_{\Q}{]}  & n \equiv 2\pmod{4} & \mbox{ is ramified at } \infty, 2 \\
	{[}(-10,1)_{\Q} {]} = {[}\Q {]} & n\in 4\Z &  \end{array} \right. 
	$$
	 whereas  
	$$
	 \Delta(H') = \left\{ \begin{array}{lll}  {[}(-10,-5)_{\Q}{]} &  n \equiv 2\pmod{4} & \mbox{ is ramified at } \infty, 5 \\
	 {[}(-10,5)_{\Q}{]}  & n\in 4\Z & \mbox{ is ramified at } 2,5 .\end{array} \right. 
	 $$
\end{example}

\subsection{Lattices, the local picture} \label{locallat}

We now assume that we work over local fields of characteristic 0. 
So let $K$ be a complete discrete valuated field with valuation ring $\Z_K$
and let $L$ be a quadratic extension of $K$ with valuation ring $\Z_L$, and let
$(V,H)$ be an $L/K$ Hermitian space. 
Denote by $\pi \in \Z_L$ a generator of the maximal ideal of $\Z_L$, if
$L/K$ is inert, then we choose $\pi \in \Z_K$ and if $L/K$ is non-dyadic ramified we choose  
$\pi $ such that $\pi^2 \in \Z_K$, i.e. $\sigma(\pi ) = -\pi $.

The next proposition is certainly well known, we indicate a proof 
as we need similar ideas for lattices invariant under a group. 

\begin{proposition} \label{duallat}
	\begin{itemize}
		\item[(a)]
	There is a $\Z_L$-lattice $\Lambda $ in $V$ 
	such that 
			$\pi \Lambda^* \subseteq \Lambda \subseteq \Lambda ^*$. Such lattices are called {\em square free}. 
\item[(b)]
For any square free lattice $\Lambda $ the Hermitian form induces a 
non-degenerate form 
$$ \overline{H} : \Lambda /\pi \Lambda^* \times \Lambda /\pi \Lambda^* \to \Z_L / \pi \Z_L, \overline{H}(\ell + \pi \Lambda ^* , m + \pi \Lambda ^*) := 
H(\ell , m) + \pi \Z_L$$
which is Hermitian in the case that $L/K$ is unramified and 
symmetric if $L/K $ is ramified. 
\item[(c)] 
If $L/K $ is ramified, then there is a $\Z_L$-lattice $\Lambda $ in 
$V$ such that $\Lambda = \Lambda ^*$. 
For such a lattice the discriminant of the symmetric bilinear space 
$(\Lambda /\pi \Lambda ^*,\overline{H})$ 
from $(b)$ is congruent to the discriminant of 
a Gram matrix of $H$ with respect to a basis of $\Lambda $.
\item[(d)]
	Let $L/K$ be inert and 
let $\Lambda $ be a square free lattice. Then the valuation 
of $\disc (H)$ is congruent modulo 2 to the dimension of the $\Z_L/\pi\Z_L$-space $\Lambda ^*/\Lambda $. 
\item[(e)]
For any  square free lattice $\Lambda $ 
$$\tilde{H} : \Lambda ^*/\Lambda \times \Lambda^* /\Lambda \to \Z_L/\pi \Z_L , \tilde{H}(\ell +  \Lambda  , m + \Lambda ) := 
       \pi H(\ell , m) + \pi \Z_L$$
	 is a Hermitian form in the case that $L/K$ is unramified and 
	a skew-symmetric form if $L/K $ is non-dyadic and ramified. 
	\end{itemize}
\end{proposition}

\begin{proof}
	(a) Let $\Lambda $ be an integral lattice, i.e. $H(\Lambda ,\Lambda ) \subseteq \Z_L$.
	Then $\Lambda \subseteq \Lambda ^*$ and $\Lambda^*/\Lambda $ is a finitely generated torsion $\Z_L$-module. 
	Let $\ell \in \Lambda ^*$ be such that $\pi ^2 \ell \in \Lambda $. 
	As $\sigma(\pi )$ is again a uniformizer of $\Z_L$ there is a unit $u\in \Z_L^{\times }$ such that 
	$$ H(\pi \ell , \pi  \ell )  =  u H ( \ell , \pi ^2 \ell ) \in \Z_L .$$
	In particular $\Lambda + \Z_L (\pi \ell ) $ is again an integral lattice. Continuing like this, we eventually 
	arrive at a lattice as in (a). 
	\\
	(b) follows from explicit computations.
	\\
	(c) We continue with the computations in (a) and 
	assume that $\ell \in \Lambda ^*$ such that $\pi \ell \in \Lambda $. 
	Then 
	$$\pi H(\ell, \ell) = H(\pi \ell, \ell )  \in \Z_L  $$
	so $\nu _{\pi }(H(\ell ,\ell)) \geq -1 $. 
	As $H(\ell ,\ell ) \in K$ the $\pi $-adic valuation 
	$\nu _{\pi }( H(\ell ,\ell)) $ is even, and therefore 
	$H(\ell, \ell) \in \Z_L$. 
	So $\Lambda + \Z_L \ell $ is again an integral lattice and 
	we can continue enlarging $\Lambda $ until $\Lambda = \Lambda ^*$.
	For a $\Z_L$-basis $B$ of such a unimodular 
	lattice $\Lambda $, the Gram matrix $H_B$ if congruent 
	mod $\pi $ to the Gram matrix of 
	the bilinear space $\Lambda/ \pi \Lambda $ and so is its determinant.
	\\
	(d) is clear and (e) follows again by direct 
	computations.
\end{proof}

\section{Algebras with involution} \label{alginvo}

I thank the referee for the motivation to set up a more sophisticated 
theory for the definition of discriminant algebras. 
%computations described in Section \ref{methods} to 
%obtain the unitary discriminants of ordinary characters.

\subsection{Central simple algebras with involution}  \label{defcsai}

This section briefly recalls some notations for algebras 
with involution over fields of characteristic $\neq 2$.
A much more sophisticated treatment 
can be found in \cite{BookOfInvolutions}.

\begin{definition}\label{definvoalg} 
	Let $L$ be a field of characteristic $\neq 2$ and 
	let $(A,\iota) $ be a central simple $L$-algebra 
	with an involution $\iota $, i.e. a ring anti-automorphism 
	of order $2$. 
	Then $\dim_L(A) = n^2$ is a square.
	\begin{itemize} 
		\item[(i)] The {\em degree} of $A$  is 
			$\deg(A) := n$. 
		\item[(ii)]  The involution is said to be of the
			first kind, if $\iota _{|L} = \id _{L} $
			and of the second kind, if the 
			restriction to the center $L$ is non-trivial.
		\item[(iii)] Put $A^+:=\{ a\in A \mid \iota(a) = a \}$ 
			and $A^{-} := \{ a\in A \mid \iota (a) = - a\} $. 
			Then $A=A^{+} \oplus A^{-}$.
		\item[(iv)] An involution of the first kind is 
			called {\em orthogonal}, if $\dim_L(A^-) = n(n-1)/2$
			and {\em symplectic}, if $\dim_L(A^-) = n(n+1)/2$.
	\end{itemize}
\end{definition}

The involutions of  matrix rings $A=L^{n\times n}$ 
are all adjoint involutions of Hermitian, symmetric or skew-symmetric forms, 
and these are of the second kind, respectively orthogonal or symplectic.  

For  any central simple $L$-algebra $A$
there is a field extension $E$, such that $E\otimes_L A \cong E^{n\times n}$.
 For $a\in A \subseteq E\otimes _L A \cong E^{n\times n}$ we define
 the {\em reduced norm} of $a$ as the determinant of the matrix
 $a\in E^{n\times n}$, $N_{red}(a) = \det(a)$ (see \cite[Section 9]{Reiner}).
 This definition depends neither on the chosen splitting field $E$ nor on the
 identification with a matrix ring.
 Moreover for any $a\in A$ we have $N_{red}(a) \in L$.

 Discriminants of involutions are only well defined if the degree of 
 $A$ is even:

\begin{remark} (\cite[Definition 7.2]{BookOfInvolutions}, \cite{Parimala})\label{discInvo}
        $A$ be a central simple $L$-algebra
        of even degree $2n$ and $\iota $ be 
	an orthogonal involution  on $A$.
        Then $A^{-}$ contains invertible elements.
	The {\em discriminant} of $\iota $ is
        defined as the square class
        $$\disc(\iota ):= (-1) ^n N_{red} (a) (L^{\times })^2 $$
        for any invertible $a\in A^{-}$.

        In the split case, $A=L^{2n\times 2n}$ and
        a symmetric $B\in \GL_{2n}(K)$ such that 
	$$\iota = \iota_B, X\mapsto B X^{tr} B^{-1}$$
        is the adjoint involution with respect to $B$,
	we have
        $$\disc(\iota) = \disc(\iota _B) = \disc (B) .$$
\end{remark}

Similar to the orthogonal case, there is a 
well defined discriminant algebra for even degree algebras with
unitary involutions:

\begin{definition}\label{discalgdef}
	(\cite[Definition (10.28), Corollary (10.35)]{BookOfInvolutions})
Let ${\mathcal A}$ be a central simple $L$-algebra of even degree $2n$
and with an involution
$\iota $ of the second kind
	Let $K:=\Fix_L(\iota _{|L})$ denote the fixed field of $\iota $ in $L$.
	Then we put the {\em discriminant algebra} 
	$[\Delta (({\mathcal A},\iota )) ]\in \Br(K)$
	of $({\mathcal A},\iota) $ 
	to denote 
	the  Brauer class of the central simple $K$-algebra  $\Delta (({\mathcal A},\iota))$
with involution 
described in 	\cite[Section 10]{BookOfInvolutions}.
\end{definition}

The discriminant algebra $\Delta (({\mathcal A},\iota )) $ 
is constructed in \cite[Section 10]{BookOfInvolutions} as a 
$K$-subalgebra of the central simple $L$-algebra 
$\lambda ^n {\mathcal A} = L \otimes _K \Delta (({\mathcal A},\iota))$. 
The latter is Brauer equivalent 
to $[{\mathcal A}]^n \in \Br(L) $. 
This algebra is split, 
$[{\mathcal A}]^n = [L] \in \Br(L) $,  if and only if $L$ splits $\Delta (({\mathcal A},\iota ))$. 
For number fields $L$ this implies that 
$L$ is a maximal subfield of 
 the quaternion algebra in $[\Delta (({\mathcal A},\iota )) ]$. 
 Then the 
 discriminant of $({\mathcal A},\iota ) $ is defined as the 
 $L$-discriminant of $[\Delta (({\mathcal A},\iota )) ]$.

 \begin{definition}\label{defdisc}
	 Assume that $L$ is a number field. 
 A central simple $L$-algebra ${\mathcal A}$ of degree $2n$ 
 is said to be {\em quasi-split}, if $[{\mathcal A}]^n = [L] \in \Br(L) $. 
	 \\
	 For an $L/K$-Hermitian involution $\iota $ 
	 on such a  quasi-split  $L$-algebra ${\mathcal A}$ we define 
	 the {\em unitary discriminant} 
	 $$\disc (({\mathcal A},\iota )) := \disc_L([\Delta (({\mathcal A},\iota ))]) = \delta N_{L/K}(L^{\times }) $$
	 if $[\Delta (({\mathcal A},\iota ))] = [(L,\delta )_K].$
 \end{definition}

	I thank J.P. Tignol for communicating to me  that 
	\cite[Corollary 2.12]{MerkurjevTignol} shows that 
	there is an extension $F$ of $K$ splitting ${\mathcal A}$ and 
	such that 
	 the map $\Br(K) \to \Br(F) , [D] \mapsto [D\otimes _K F ]$ 
is injective. 
	This observation is also used 
	in the proof of \cite[Corollary (10.36)]{BookOfInvolutions}. 
	It allows to conclude the following remark.

\begin{remark} \label{discalggen} 
	The discriminant algebra $[\Delta (({\mathcal A},\iota)) ]$ is the unique element
	of $\Br(K)$ with the following property: 
 For any extension $E^+$ of $K$ such that  $E:=E^+L $ is 
a splitting field of ${\mathcal A}$ and a Hermitian 
form $H$ on $E^{2n}$ with $(\iota _H)_{|{\mathcal A}} = \iota $  we get 
$$[E^+\otimes _K \Delta (({\mathcal A},\iota ))] =  
	= \Delta (H)  = [ (E,\disc(H))_{E^+} ] \in \Br(E^+) .$$
\end{remark}

 Note that  Remark \ref{discalggen} characterises the 
 unitary discriminant of a quasi-split central simple $L$-algebra
 $({\mathcal A},\iota )$ by the following property: 
	For 
	any extension $E^+$ of $K$ such that  $E:=E^+L $ is 
	a splitting field of ${\mathcal A}$ and a Hermitian 
	form $H$ on $E^{2n}$ with $(\iota _H)_{|{\mathcal A}} = \iota $  we get 
	${\mathcal A}\otimes E = E^{2n\times 2n}$ and 
	 $$ \disc(({\mathcal A},\iota ) )  N_{E/E^+}(E^{\times }) = \disc(H) .$$

 %\begin{proof}
	 %Let $\delta \in K^{\times }$ represent $\disc(({\mathcal A},\iota ) ) $.
	 %By Remark \ref{discalggen}
	 %we get $$
	 %[(E,\delta )_{E^+} ] = 
	 %[\Delta ({\mathcal A},\iota ) \otimes _{K} E^+] = 
	 %[\Delta (H) ] = [(E,\disc(H) )_{E^+}]  $$
	 %so $\delta N_{E/E+}(E^{\times }) = \disc(H) $.
 %\end{proof}

\subsection{Orthogonal subalgebras}\label{ind+}

\begin{definition}\label{deforthsub}
Let ${\mathcal A}$ be a central simple $L$-algebra of even degree
and with $L/K$-involution
$\iota $ of  second kind.
        Then a $K$-subalgebra $A$ of ${\mathcal A}$ is called an
        {\em orthogonal  subalgebra} of $({\mathcal A},\iota )$,
        if and only if
        \begin{itemize}
                \item[(a)] $A$  is a
                        central simple $K$-algebra with $LA = {\mathcal A}$.
                \item[(b)] $A$ is invariant under $\iota $, i.e. $\iota(A) = A$.
                \item[(c)] The restriction of $\iota $ to $A$
                        is an orthogonal involution of $A$.
        \end{itemize}
\end{definition}

Note that orthogonal subalgebras do not exist in general.
However for split algebras, there are always orthogonal subalgebras.

\begin{example}\label{orthsubexample}
        Let $\iota $ be the adjoint involution of a Hermitian form $H$
        on $V$
        and ${\mathcal A}:=\End_L(V)\cong L^{n\times n}$.
        Then for any orthogonal basis $B$ of $H$
        the subalgebra
        $$A:=\{ a\in {\mathcal A} \mid \ ^B a ^B \in K^{n\times n } \} \cong
        K^{n\times n} $$
        is an orthogonal subalgebra of $({\mathcal A},\iota )$.
\end{example}

\begin{theorem}\label{orthsub} (\cite[Proposition (10.33)]{BookOfInvolutions}, \cite[Proposition 11]{Queg})
Let ${\mathcal A}$ be a central simple $L$-algebra of even degree $2n$
and with $L/K$-involution
$\iota $ of  second kind.
	\begin{itemize}
		\item[(a)]
	Let $A$ be an orthogonal subalgebra of $({\mathcal A},\iota )$.
	Then 
			$$[\Delta(({\mathcal A},\iota ))] = [A]^n [(L,\disc(\iota _{|A }))] \in \Br(K).$$
\item[(b)] 
	If $A$ is a symplectic subalgebra of $({\mathcal A},\iota )$, then 
			$$[\Delta(({\mathcal A},\iota ))] = [A]^n  .$$
	\end{itemize}
\end{theorem} 

If $({\mathcal A},\iota)  \cong (L^{2n \times 2n},\iota _H)$ is split, then the formulas in 
Theorem \ref{orthsub} applied to the discriminant $\disc(H)$ of the 
Hermitian form read as 
\begin{itemize}
	\item[(a)] $\disc(H) =(-1)^n \disc_L([A]^n) N_{red}(X) $ 
		for any $X=-\iota(X) \in A^{\times }$ in the orthogonal case and 
	\item[(b)] $\disc(H) = \disc_L([A]^n)$ in the symplectic case. 
\end{itemize} 

\subsection{Stable subalgebras} \label{semisimple}

We focus on involutions of the second kind, similar results also hold for 
orthogonal involutions (see \cite[Theorem 5.13]{OrthogonalStability}). 
        So let $L$ be a field of characteristic $\neq 2$ and let
	${\mathcal A}$ be a central simple $L$-algebra
        of even degree and with an involution $\iota $ of the second kind.

	\begin{definition}
		A semisimple subalgebra $A\leq {\mathcal A}$ is called {\em stable} 
		if 
		\begin{itemize} 
			\item $\iota (A) = A$ and
			\item all $\iota $-invariant simple direct summands of 
				$A$ have even degree. 
		\end{itemize} 
	\end{definition}

So a semisimple involution invariant subalgebra $$ A = \bigoplus _{i=1}^h A_i = 
 \bigoplus _{i=1}^h e_i A  $$ 
 is stable, if for all central primitive idempotents $e_i$ with $\iota (e_i) = e_i$, 
 the dimension $\dim_{Z(A_i)} (A_i)$ is even. 
If $\iota (e_i ) = e_j$ for some $j\neq i$ 
then the restriction of $\iota $ to $(e_i+e_j)A$ is hyperbolic and has 
no contribution to the discriminant algebra of $({\mathcal A},\iota)$. 
 The next remark shows that, up to hyperbolic involutions, 
 we can reduce to the simple case by replacing ${\mathcal A}$ by 
 $e_i{\mathcal A} e_i$. 

\begin{remark} \label{simple}
        If ${\mathcal A} = {\mathcal D}^{n\times n}$ is a central simple $L$-algebra
        (${\mathcal D} $ a central $L$-division algebra) and $e^2=e\in {\mathcal A}$,
        then $e{\mathcal A} e \cong {\mathcal D}^{x \times x}$ for some $x\leq n$.
        An easy proof uses the fact that ${\mathcal A}$ has a unique simple module,
        ${\mathcal D}^n$ and $e{\mathcal A}$ is a direct sum of simple modules, so
        isomorphic to $\oplus ^x {\mathcal D}^n$. Then
        $$(e{\mathcal A} e)^{op} \cong \End_{\mathcal A} (e {\mathcal A} )
        \cong ({\mathcal D}^{op})^{x\times x} $$
        implies that $e{\mathcal A} e \cong {\mathcal D}^{x \times x}$.
\\
        So if $e^2=e=\iota(e) \in A$ is an involution invariant central
        primitive idempotent 
        \begin{itemize}
                \item[(a)] $eAe$ is a simple algebra in  $e{\mathcal A}e  $ and
                \item[(b)] $e{\mathcal A}e \cong {\mathcal D}^{x\times x}
                        $ is a central simple $L$-algebra Brauer equivalent to
                        ${\mathcal A}$.
        \end{itemize}
\end{remark}

To handle reducible characters, we need to define the discriminant
of a module. 
The idea is that the discriminant algebra of a central simple 
stable algebra with involution is the discriminant algebra of its unique simple 
module.  
Further this discriminant algebra behaves multiplicative with respect to 
direct sums of modules and hence it is enough to deal with the case that
$A$ is a simple stable subalgebra of ${\mathcal A}$.

\begin{remark}\label{simplesub}
	Assume that $A$ is a simple stable subalgebra of ${\mathcal A}$
	such that $Z(A) \leq Z({\mathcal A}) = L$. 
	Then $LA$ is a central simple $L$-algebra with an involution  $\iota _{|LA}$
	of the second kind. 
	Let $m$ be the multiplicity of the simple $LA$-module in 
	the restriction of the simple ${\mathcal A}$-module to $LA$. 
	%Let $C:=C_{{\mathcal A}}(LA)  $ be the centraliser algebra 
	%of $LA$ in ${\mathcal A}$. 
	%Then there are orthogonal primitive $\iota $-invariant 
	%idempotents $e_1,\ldots , e_m$ in $C$ with $1=e_1+\ldots + e_m$. 
	%For all $1\leq i \leq m$ the algebra 
	%$e_i LA e_i \leq e_i {\mathcal A} e_i$
	Then $$[\Delta(({\mathcal A},\iota ))] = 
	[\Delta ((LA, \iota _{|LA }))] ^m .$$ 
\end{remark}

It remains to consider the situation where $E:=Z(A) > Z({\mathcal A})=L$. 
Here we assume that the characteristic of $L$ is $0$, to avoid inseparable 
extensions. 
Let $F:=\Fix_{\iota }(E)$ denote the fixed field of $\iota $ in $E$.
Then $F\geq K$ is a separable finite extension of $K$ and hence 
there is a well defined map, known as corestriction \cite{Tignol}, 
$$ N_{F/K} : \Br (F) \to \Br (K)  .$$
In particular $N_{F/K} ([(E,\delta )_F]) = [(L,N_{F/K}(\delta ))_K]$ 
by \cite[Theorem 3.2]{Tignol} (and the fact that $E=LF=F[\ell ]$ for some $\ell \in L$).

\begin{proposition}\label{correst}
	Let $({\mathcal A},\iota ) $ be a central simple $L$-algebra
        of even degree  and with $L/K$-Hermitian involution $\iota $.
	Let $L\leq E \leq  {\mathcal A}$ be an extension field of $L$,
        such that $\iota(E) = E$. Let $F$ denote the fixed field of $\iota $ in $E$.
	Put $ A:= C_{{\mathcal A}}(E)$ to denote the centraliser of 
	$E$ in ${\mathcal A}$.
        Then $(A,\iota ) $ is a central simple $E$-algebra with involution.
        If the degree of $A$ is even, then
	$$[\Delta (( {\mathcal A},\iota ) )]  = N_{F/K} ([\Delta( (A,\iota ))]) .$$
\end{proposition}

\begin{proof}
	 By Remark \ref{discalggen} we can assume  that 
 ${\mathcal A} \cong L^{2n \times 2n} $ is split.
	 Then also $A \cong  E^{2d \times 2d}$ is split, where $n/d = [E:L]$. 
	 So there is an $E/F$-Hermitian form $H\in A$ such that $\iota _{|A} = \iota _H$. 
	 Let  $\varphi : E\hookrightarrow L^{n/d\times n/d}$ denote 
	 the regular representation. 
	 Applying $\varphi $ to the entries of the matrices yields an embedding 
	 $\phi $ of $E^{2d \times 2d} $ into $L^{2n\times 2n} $. 
	 Then $\det(\phi (H)) = N_{F/K} (\det(H)) $, so 
	 $$[\Delta (( {\mathcal A},\iota ) )] =  
	 [(L,\disc(\phi(H)) _K] = [(L,N_{F/K}(\disc(H)))_K] = N_{F/K} ([\Delta( (A,\iota ))])  .$$ 
\end{proof}

\begin{corollary}\label{lemnorm}
	Assume that $A$ is a simple stable subalgebra of ${\mathcal A}$
	such that $E = Z(A) \geq Z({\mathcal A}) = L$ and let $F:=\Fix_{\iota }(E)$.
	Assume that $E/K$ and  $F/K$ are Galois extensions.
	Let $m$ be the multiplicity of the simple $A$-module in 
	the restriction of the simple ${\mathcal A}$-module to $LA$. 
	Then $m = [F:K] m'$ is a multiple of $[F:K]$ and 
	 $$[\Delta(({\mathcal A},\iota ))] = 
	  N_{F/K}([\Delta ((A,\iota _{|A }))]) ^{m'} .$$
	 In particular, if ${\mathcal A}$ and $A$ are quasi-split then 
	 $$\disc(({\mathcal A},\iota )) = N_{F/K} (\disc((A,\iota )))^{m'} $$
	 where $N_{F/K}$ is the usual norm map.
\end{corollary}

%\begin{proof}
	%As 2 is invertible there is a system $e_1,\ldots ,e_{m'}$ of 
	%involution invariant idempotents in $C_{{\mathcal A}}(A)$. 
	%Replacing ${\mathcal A}$ by $e_i {\mathcal A} e_i$ and 
	%$A$ by $e_i A e_i$ we are hence in the situation of Proposition 
	%\ref{corest}. 
%\end{proof}

\section{Discriminants of characters}  

\subsection{Characters of finite groups} 

Let $G$ be a finite group and $K$ be a field. 
A $K$-representation $\rho $ of $G$ is a group homomorphism 
$\rho : G \to \GL_n(K)$ into the group of invertible $n\times n $-matrices
over $K$. 
The associated Frobenius character is the map 
$$\chi_{\rho }: G \to K, g\mapsto \trace (g) .$$
For fields $K$ of characteristic zero, Frobenius characters 
are in bijection to isomorphism classes of $K$-representations of $G$. 
As $\rho(g)$ is diagonalizable over an algebraic closure of $K$, 
the character values are sums of $|G|$-th roots of unity. 
In particular the {\em character field} $\Q(\chi_{\rho})$ of $\chi _{\rho }$,
the field extension $\Q(\chi_{\rho }(g) : g\in G )$ generated by
the character values, is an abelian number field. 
A character is called {\em irreducible} if it is not a proper sum of 
characters. 
The {\em ordinary character table} 
displays the character values of the complex irreducible characters 
of $G$ on representatives of the conjugacy classes of $G$. 
For more details see for instance \cite{Isaacs}.

To relate the representation theory over fields of characteristic 
$0$ and positive characteristic, $p$, Brauer characters are introduced. 
An element of $G$ is called a $p'$-element, if its order is not 
divisible by $p$. The set $G_{p'}$ denotes the set of all $p'$-elements
in $G$; it is a union of conjugacy classes, the $p'$-classes.
A $p$-Brauer character is a certain class function $\varphi : G_{p'} \to \C $.
Loosely speaking,
the restriction $\chi \pmod{p} $ 
of a complex Frobenius character $\chi $ to the $p'$-classes of
$G$ is a $p$-Brauer character of $G$. 
If $p$ does not divide the group order, then 
restriction defines a bijection between the 
irreducible $p$-Brauer characters and  the irreducible complex characters. 
For primes $p$ dividing the group order the 
$p$-Brauer character $\chi \pmod{p}$ might be reducible.
The irreducible constituents of 
$\chi  \pmod{p} $ are in bijection to the 
composition factors of $L/\wp L$, for any lattice 
$L$ invariant under the representation $\rho $ affording $\chi $ 
and a suitably chosen prime ideal $\wp $. 
To specify the ideal we also denote this modular reduction
by $\chi \pmod{\wp}$.
See \cite[Sections 6 and 7.1]{OrthogonalStability} 
and \cite[Introduction]{BrauerATLAS}
for the subtleties that might arise here and 
Section \ref{ONan} for an example.
For a more sophisticated introduction to Brauer characters 
see \cite[Chapter 15]{Isaacs} or, for a summary,
\cite[Section 2.1]{survey} and \cite[Introduction]{BrauerATLAS}.
We just note the following.

\begin{remark}\label{Brauerclassification} 
Let $F$ be a sufficiently large 
	finite field of characteristic $p$.
	Then the irreducible $p$-Brauer characters are in bijection
	to the simple $FG$-modules.
	We then also say that the corresponding $F$-representation of 
	$G$ affords the respective irreducible Brauer character.
\end{remark}

\begin{definition} (see for instance \cite[Introduction, Section 9]{BrauerATLAS})
	Let $\chi $ be an
irreducible ordinary 
	or  Brauer character of the finite group $G$ and let 
	$\rho $ be a representation affording $\chi $.
	The Frobenius-Schur indicator of $\chi $ is
	\begin{itemize} 
		\item[+] $\ind(\chi ) = +$ if $\chi $ is real valued and 
			there is a non-degenerate 
			$\rho(G)$-invariant quadratic form 
			or $\chi $ is the trivial 
			$2$-Brauer character. 
		\item[-] $\ind(\chi ) = -$ if $\chi $ is real valued and 
			there is a non-degenerate 
			$\rho(G)$-invariant symplectic form 
			but no invariant non-degenerate quadratic form. 
		\item[o] $\ind(\chi ) = o$ in all other cases.
	\end{itemize}
\end{definition}

Note that the set of invariant quadratic, symplectic, or Hermitian 
forms is a vector space. In particular also scalar multiples 
of invariant forms are invariant. 
As $\det(a H) = a^{\chi(1)} \disc( H) $ the discriminant of the
invariant quadratic or Hermitian forms is only well-defined if 
the degree of the character $\chi $ is even: 

\begin{definition}
	For $s\in \{ +, o \} $ the set 
        $$\Irr ^{s}(G) := \{ \chi \in \Irr_{\C }(G) \mid
        {\mathrm{ind}}(\chi ) = s, \chi(1) \in 2\Z \} $$
denotes the even degree irreducible ordinary characters of
$G$ which have Frobenius-Schur indicator $o$ respectively $+$.
\end{definition}

The Frobenius-Schur indicator of an ordinary irreducible character $\chi $ is
easily calculated from its character values. 
If the Frobenius-Schur indicator of $\chi $ is $o$  then 
there is a positive-definite $\rho(G)$-invariant Hermitian form 
for any representation $\rho $ affording $\chi $.
This is not automatically true for Brauer characters,
as there are irreducible Brauer characters of indicator $o$ for which 
the associated representation is not contained in a unitary group, 
for instance all indicator $o$ irreducible 3-Brauer characters 
of $L_3(3)$ are not unitary, as  all irreducible $3$-modular
representations of $L_3(3)$ are realisable over $\F_3$. 

\begin{definition}
	An irreducible $p$-Brauer character $\chi $ is called 
	{\em unitary}, if there is a $p$-power $q$ such that 
	$\chi (g^{-1}) = \chi(g^{q} ) $ for all $g\in G_{p'}$. 
\end{definition}

Note that real valued Brauer characters 
(indicator $+$ or $-$) are always unitary as 
$\chi(g^{-1}) = \overline{\chi}(g)$ so $q=p^0=1$ is 
a possible choice. 
The unifying philosophy here is to consider 
symplectic and quadratic forms as ``unitary'' invariant forms over a 
suitable field extension. Any such non-degenerate form gives rise to 
an isomorphism between the representation and its 
(unitary) dual.

\begin{lemma} (see \cite[Lemma 4.4.1]{Colva})
	An indicator $o$  irreducible 
	Brauer character $\chi $ is unitary if and only if
	the representation $\rho $ affording $\chi $ 
	admits a non-degenerate unitary form.
\end{lemma} 

%\begin{proof}
%Let $\rho : G\to \GL_n(\F_{q^2})$ be a representation 
%affording the character $\chi $ and denote by $\phi $ the Frobenius 
	%automorphism of $\F_{q^2}$ with fixed field $\F_q$.
%Then any non-degenerate $\rho(G)$-invariant unitary form $H$ satisfies 
%$$\rho(g) H \phi (\rho(g))^{tr} = H \mbox{ for all } g\in G $$
%hence $H$ conjugates 
	%$\rho(g)$ to $\hat{\rho}(g):=\phi (\rho(g^{-tr}) $. 
%In particular these two simple modules are isomorphic and hence
%have the same Brauer character.
%For any $p'$-element $g\in G$ the trace of $\phi (\rho(g))$
	%is the sum of the $q$-th powers of the eigenvalues of $\rho (g)$,
	%i.e. the trace of $\rho(g^q)$. 
	%So the Brauer character $\hat{\rho }$ is the map
	%$g\mapsto \chi(g^{-q}) $.
	%\\
	%For the other direction 
%let $q$ be a $p$-power such that 
	%$\chi(g^{-1}) = \chi(g^q)$ for all $g\in G_{p'}$. 
	%Then $\chi(g) = \chi(g^{q^2})$ 
	%for all $g\in G_{p'}$, so the $\F_p(\chi ) \subseteq \F_{q^2}$.
	%Let $\rho: G\to \GL_n(\F_{q^{2}})$ 
	%be the representation affording $\chi $.
	%Define 
	%$$\hat{\rho}:G\to \GL_n(\F_{q^{2}}) , \hat{\rho}(g) := \phi (\rho(g^{-tr})  $$ 
	%where $\phi $ if the Frobenius automorphism of $\F_{q^2}/\F_q$. 
	%Then these two representations have the same irreducible 
	%Brauer character, so they are equivalent, i.e. there is 
	%$H\in \GL_n(\F_{q^2})$ with 
	%$$\rho(g) H \phi (\rho(g))^{tr} = H \mbox{ for all } g\in G .$$
%Then standard arguments show that 
	%$\phi(H)^{tr} H^{-1} $ commutes with $\rho(G)$ and hence is 
	%a scalar matrix, so $H$ is equivalent to some unitary form.
%\end{proof}

	To simplify the exposition we also say that  
	all ordinary characters are unitary.

\subsection{Schur indices} 

Schur indices play an important role in the computation of unitary discriminants 
of characters. As described in Theorem \ref{orthsub} we need all 
local Schur indices to compute these discriminants.
These are encoded in the Brauer element, a notion coined in 
\cite[Definition 2.1]{Turull}:

\begin{definition}
	Let $\chi $ be an absolutely irreducible ordinary character 
	of some finite group $G$. 
	Let $K$ be some field containing the 
	character field ${\mathbb{Q}}(\chi )$  and 
         let $\rho $ be a 
	$K$-representation of $G$ affording the character  $m\chi $ for some positive integer $m$.
	Then $m$ is minimal if and only if $\End_{KG}(\rho )  =: D$ is 
	a central simple division algebra over $K$. 
	In this case $m^2 = \dim_K(D)$ and $m_K(\chi):=m $ 
	is the {\em Schur index} of $\chi $ over $K$.
	The class $[\chi ]_K:=[D] \in \Br(K)$ of $D$ in the Brauer group $\Br(K)$ of $K$ 
	is called the {\em Brauer element} of $\chi $ over $K$.
	\\
	The field $K$ is called a {\em splitting field} of $\chi $ 
	if the Brauer element $[\chi ]_K = [K]$ is  trivial.
	\\
	If $K=\Q(\chi )$ then we often omit the index $K$ and 
	talk about the Schur index of $\chi $. 
\end{definition}

It is well known that the Schur index $m$ of an irreducible character $\chi $ 
always divides the degree $\chi (1)$. 

\begin{definition} 
	A character $\chi \in \Irr^o(G)$ with Schur index $m$ is called 
	{\em quasi-split}, if 
 $\chi(1)/m$ is even. 
\end{definition}

\subsection{The unitary discriminant of a character} 

The Galois orbits of irreducible ordinary characters of a finite group $G$ 
are in bijection to the direct summands of the rational group algebra 
$\Q G$. More precisely let $\chi_1,\ldots , \chi _h \in \Irr (G)$ represent these
orbits then 
$$\Q G = \bigoplus _{i=1}^h A_i $$ 
where the center $L_i:=Z(A_i) = \Q(\chi _i)$ is the character field of $\chi _i$ 
and $\chi_i(1)$ is the degree of the central simple $L_i$-algebra $A_i$, i.e. 
$\chi _i(1)^2 = \dim _{L_i}(A_i) $. 
The natural involution $\iota $ on $\Q G$, defined by $\iota (g) = g^{-1}$ for
all $g\in G$, preserves the simple summands of $\Q G$. 
For $1\leq i \leq h$ let $\iota _i := \iota _{|A_i }$ denote the restriction of $\iota $ to $A_i$.

\begin{remark}
For $1\leq i \leq h$ the involution $\iota _i$ is 
	\begin{itemize} 
		\item symplectic, if  $\ind (\chi _i) = - $, 
		\item orthogonal, if $\ind (\chi _i) = +$,  and 
		\item of the second kind (or unitary), if  $\ind(\chi _i) = o$.
	\end{itemize}
\end{remark}

\begin{definition}\label{defdiscrchi}
	For $\chi _i \in \Irr ^+ (G) $ 
	$$\disc(\chi _i):= \disc((A_i,\iota _i)) \in L_i^{\times }/(L_i^{\times })^2 $$ is called 
	the {\em orthogonal discriminant} of the character $\chi _i$.
	\\
	For $\chi _i \in \Irr ^o (G) $ the character field 
	$L_i $ is a complex number field and the 
	restriction of $\iota _i$ to $L_i$ is complex conjugation. 
Let $K_i:=L_i^+$ the maximal real subfield of $L_i$.
	Then
	$$\Delta (\chi _i) := [\Delta((A_i,\iota _i))]  \in \Br(K_i)$$ is called 
the {\em unitary discriminant algebra} of the character $\chi _i$.
\\
	If $\chi _i \in \Irr^o(G)$ is quasi-split, then 
	$A_i$ is a quasi-split central simple $L_i$-algebra.  Then we call 
	$$\disc(\chi _i ):= \disc ((A_i,\iota _i)) \in K_i^{\times }/N_{L_i/K_i}(L_i^{\times }) $$ 
	the {\em unitary discriminant} of $\chi _i$.
\end{definition}

The discriminant algebra of a unitary involution is usually hard to 
compute. It is however possible for quaternion algebras, as explained 
in \cite[p. 129]{BookOfInvolutions}:

Consider the degree 2 faithful irreducible character $\chi $ of 
$C_7 \times Q_8$. Then $L=\Q(\chi ) = \Q[\zeta_7]$ is 
the 7th cyclotomic number field and the Brauer element 
$[\chi ]_L = [(-1,-1)_L]$ is the class of the quaternion algebra over $L$ 
ramified at the two prime ideals dividing 2.
By \cite[(2.22)]{BookOfInvolutions} the discriminant algebra of 
this character is $[(-1,-1)_{\Q(\zeta _7 + \zeta _7)} ]$.
Note that here $\chi $ is not quasi-split, so it does not have
a unitary discriminant.

Also if the Schur index of $\chi \in \Irr^o(G)$ is odd, 
then the discriminant algebra of $\chi $ can be obtained by computing 
invariant Hermitian forms after passing to an odd degree 
extension of the character field. 
If two Hermitian spaces become isometric over an odd degree extension, 
then   \cite[Corollary 6.16]{BookOfInvolutions} tells us that they are 
already isometric over the base field.

For instance consider the famous example of the faithful degree 
$3$ character $\psi $ of $H:=C_7  \rtimes C_9$ with center $C_3$ and 
character field $\Q(\psi ) = L:=\Q[\sqrt{-7}, \sqrt{-3} ]$.
Then $\psi $ has  Schur index $3$.
We tensor $\psi $ 
with the degree 2 character of the group $D_8$ to obtain a 
degree $6$ irreducible character $\chi \in \Irr^o(H\times D_8)$ with
character field $\Q(\chi ) = L$. 
There is a monomial representation of $H\times D_8$ over 
$E:=\Q(\zeta_{21} )$ affording the character $\chi $. 
Over this degree 3 field extension the unitary discriminant of 
$\chi $ is $-1$, so we find 
$$\Delta (\chi ) \otimes _{L^+} E^+ = [(E, -1 )_{E^+} ] $$
and hence $\Delta (\chi ) = [(L, -1)_{L^+}]  $ is the 
quaternion algebra over $L^+ = \Q(\sqrt{21})$ ramified only at the 
two infinite places of $L^{+}$. 
As $\chi $ is quasi-split, it has a well defined unitary discriminant,
$\disc(\chi ) = -1$.

\subsection{Unitary stability}

In \cite{OrthogonalStability} the most important notion is the 
one of orthogonal stability. 
 An ordinary or Brauer character $\chi $ is called 
 {\em orthogonally stable}, if all its indicator + constituents 
 have even degree. 
 Orthogonally stable characters are exactly those that have
 a well defined orthogonal discriminant. 
For unitary characters we need a slightly stronger condition:

\begin{definition}\label{unitarystable}
	An ordinary or Brauer character $\chi $ of a finite group $G$ 
	is called {\em unitary stable}, if 
	all irreducible  unitary 
	constituents of $\chi $ have even degree. 
\end{definition}

It is clear that unitary stability is a necessary condition for an ordinary
character to have a well defined discriminant algebra, as the multiplication 
of the invariant forms in odd degree constituents changes the discriminant.

The idea behind the definition of unitary stability is twofold: 
On the one hand
unitary stable modular reductions allow to exclude prime divisors of the 
discriminant (see Theorem \ref{duallatconst}).
On the other hand Proposition \ref{rechus} shows that 
unitary stable ordinary characters are exactly those
that have a well defined discriminant algebra. 

Whereas orthogonal discriminants of orthogonally stable characters 
are essentially independent of the chosen splitting field, 
for unitary discriminants this field matters. 

\begin{definition}    %TODO
Let $\chi \in \Irr^o(G)$
	and let $L$ be a totally 
	complex number field containing $\Q(\chi )$ and with real subfield, say, $K$.
	The {\em $L$-discriminant algebra} of $\chi $
	is $$[\Delta _L(\chi )]:= [K\otimes _{\Q(\chi)^+} \Delta (\chi)] .$$
\end{definition}

For real characters the $L$-discriminant algebra is given in Proposition \ref{ind+-}
in the next section.

\subsection{Unitary discriminants of real characters} \label{indicator}

Let $\chi $ be an irreducible real character of $G$ of even degree 
$2n=\chi(1)$. 
Then its character field $K:=\Q(\chi )$ is a real number field, 
the Frobenius-Schur indicator is
$\ind (\chi ) \in \{+,- \}$ and the Schur index of $\chi $ is 
$m \in \{1,2 \}$. 
Let $\rho : G \to \GL_{\chi(1) m} (K) $ be a $K$-representation affording 
the character $m\chi $ and put $A := \langle \rho (g) \mid g\in G \rangle _K$ 
to denote the central simple $K$-algebra generated by the matrices in $\rho(G)$.
Then $A$ is invariant under the natural involution $\iota $ inverting the
group elements. 

Let $L = K[\sqrt{-\delta }] $ be a totally complex quadratic extension 
of $K$.
Extend the involution $\iota $ to the $L/K$-Hermitian involution 
$\widetilde{\iota }:= \gamma \otimes \iota $ on ${\mathcal A}$ where 
$\langle \gamma \rangle = \Gal (L/K)$.
Then Theorem \ref{orthsub} allows to conclude the 
following proposition.

\begin{proposition} \label{ind+-} 
	\begin{itemize} 
	\item[a)] 
	If the indicator of $\chi $ is $+$ then 
	$\chi $ is orthogonally stable and hence has 
	a well defined orthogonal discriminant $\disc(\chi )$. 
	Then 
		$$\Delta _L(\chi ) :=[\Delta ({\mathcal A},\widetilde{\iota})]
			= [(-\delta , \disc(\chi ))_K ][\chi]_K^n .$$
			%If ${\mathcal A}$ is quasisplit then 
			%$\disc _L(\chi) = \disc(\chi ) \disc_L([\chi ]_K^n) $.
	\item[b)]  
	If the indicator of $\chi $ is $-$ then 
		$$\Delta _L(\chi ) :=[\Delta ({\mathcal A},\widetilde{\iota})]
			= [\chi]_K^n .$$
			%If ${\mathcal A}$ is quasisplit then 
			%$\disc _L(\chi) = \disc_L([\chi ]_K^n) $.
		%In particular if $n$ is even then $\disc_L(\chi )$ is 
		%trivial.
	\end{itemize} 
\end{proposition}

%\begin{proof}
	%Part a) follows directly from Theorem \ref{orthsub} from which we can 
	%also deduce part b) by replacing the symplectic involution on ${\mathcal Q}$ 
	%by some orthogonal involution: 
	%Write ${\mathcal Q} := \langle 1,i,j,ij \rangle _K $ with $i^2=a$, $j^2=b$, 
	%$ij=-ji$. Embedding ${\mathcal Q}$ diagonally 
	%into $A= {\mathcal Q}^{n \times n}$ 
	%we can change the symplectic involution $\iota $ of $A$ to $\kappa _i \circ \iota $, where $\kappa _i $ is conjugation by $i$. 
	%Taking $\ell \in L$ such that $\gamma (\ell ) = -\ell $ the involution 
	%$$J :=( \kappa _{\ell } \otimes \kappa _i ) \circ \widetilde{\iota } $$ 
	%is an involution of the second kind on ${\mathcal A}$ restricting 
	%to an orthogonal involution on $A$. 
	%So $A$ is an orthogonal subalgebra of $({\mathcal A},J)$ and hence 
	%Theorem \ref{orthsub} yields that 
	%$$[\Delta (({\mathcal A},J))] = [{\mathcal Q}\otimes K]^n 
	%[(L,\disc(J _{|A }))_K] = 
	%[{\mathcal Q}\otimes K]^n  [(L,N_{red}(i))_K] $$
	%where we use the fact that $[A] = [{\mathcal Q}] \in \Br(K)$ 
	%and that the discriminant of an orthogonal involution is the reduced norm
	%of a skew symmetric element.
	%Now $\ell \otimes i $ is a symmetric element in $({\mathcal A},\widetilde{\iota})$, i.e. $\widetilde{\iota }(\ell \otimes i ) = \ell \otimes i $, \cite[Corollary (10.36)]{BookOfInvolutions} says 
%that 
	%$$
	%[\Delta(({\mathcal A},\widetilde{\iota }))] =
	%[\Delta (({\mathcal A},J))] [(L,N_{red}(\ell \otimes i))_K] = [{\mathcal Q}\otimes K]^n  $$  
%because $[(L,N_{red}(\ell \otimes i))_K] = [(L,N_{red}(i))_K] $. 
%\end{proof}

\subsection{The discriminant of a unitary stable character}

\begin{proposition}\label{rechus}
	Let $G$ be a finite group and let
	$$\Chi := \sum _{i=1}^p a_i \chi_i + \sum_{j=1}^s b_j \psi_j + 
	\sum _{k=1}^u c_k \varphi _k $$
	be an ordinary unitary stable character and written as 
	sum of irreducibles such that $\chi _i \in \Irr^+(G)$ occur 
	with multiplicity $a_i$, $\ind(\psi_j) = -$ and $\varphi _k \in \Irr^{o}(G)$. 
	$L = \Q(\Chi )$ be  the character field of $\Chi $ and  denote the maximal real 
	subfield of $L$ by $K$. 
	Assume that $L\neq K$. 
	Then there is a unique class $\Delta (\Chi ) \in \Br (K) $ 
	such that for any extension $E$ of $L$ and any
	$E$-representation $\rho : G \to \GL_{\Chi(1)}(E)$ affording the
	character $\Chi $ 
	all non-degenerate $\rho (G) $-invariant 
Hermitian forms have discriminant algebra $\Delta(\Chi )\otimes E^+ \in \Br(E^+) $.
\end{proposition}

Note that the assumption that $L\neq K$ is made to simplify the statement,
the case that $L=K$ has been dealt with in \cite[Section 5]{OrthogonalStability}.

\begin{proof}
	By Remark \ref{simple} the discriminant algebra of $\chi $ is 
	obtained as the product of the discriminant algebras 
	of the sums of the Galois orbits. 
	In particular 
	$$\Delta (\Chi ) = \Delta _L (\chi ) \Delta _L (\psi ) \Delta _L (\varphi )$$
	where 
	$\chi  = \sum_{i=1}^p a_i \chi_i $, 
	$\psi = \sum_{j=1}^s b_j \psi_j $, and $\varphi = \sum _{k=1}^u c_k \varphi _k $.
	\begin{itemize}
		\item[(a)] 
			The character $\chi $ is orthogonal and 
			orthogonally stable. By \cite[Theorem 5.13]{OrthogonalStability} 
			the character has a well defined orthogonal discriminant 
			$\disc (\chi ) $.  
			By Theorem \ref{orthsub} 
		the contribution of $\chi$ to $\Delta (\Chi )$ 
			is $$\Delta _L(\chi ) = [\chi  ]_K^{\chi(1)/2} [(L,\disc(\chi ))_K] .$$
		\item[(b)]
			For $1\leq j\leq s $ let 
			 $\Gamma_j := \Gal (L(\psi_j) /L )$ and choose a system 
			 of representatives $\psi _j , j \in J$ 
			 of the Galois orbits on the set 
			 $$\{\psi _1,\ldots , \psi _s \} = \bigcup_{j\in J} \psi_j^{\Gamma _j} .$$ 
			 Then by Theorem \ref{ind+-} b) and Corollary \ref{lemnorm} 
			 the discriminant algebra of the Galois orbit is
			 $$\Delta _j := \Delta _L(\sum_{\gamma \in \Gamma_j } \psi_j^{\gamma } )  
			 = N_{L(\psi_j)/L} ([\psi _j]^{\psi_j(1)/2}) .$$
		So the contribution of $\psi $ to $\Delta (\Chi )$  is
			$$\Delta _L(\psi ) = \prod_{j\in J} \Delta _j^{b_j}.$$
			 \item[(c)] 
				 As in (b), we choose a system of representatives
			$\{\varphi _k \mid k \in \Kappa \} $ of the 
			Galois orbits and put for $k\in \Kappa $ 
			$$\Delta _k := \Delta _L(\sum_{\gamma \in \Gamma_k} \varphi_k^{\gamma } ) = N_{L(\varphi _k)/L} (\Delta_{L(\varphi_k)} (\varphi _k) ) $$ 
			and use Corollary \ref{lemnorm} to obtain 
			$$\Delta _L(\varphi )  =  \prod _{k\in \Kappa } \Delta _k^{c_k} .$$
	\end{itemize}
\end{proof}

\section{Elementary character theoretic methods} \label{methods}

In this section we let $G$ be a finite group, 
$\chi \in \Irr^o(G)$, $L=\Q(\chi )$, and $K=L^+$.  
This section describes elementary character theoretic methods that 
are used to  compute the unitary discriminant algebra of $\chi $.

 Remark \ref{infty} 
	translates into the following remark.

\begin{remark}\label{infinitechi}
An infinite place of $K$ ramifies in 
	$\Delta (\chi )$
	if and only if $\chi (1) \equiv 2 \pmod{4}$. 
\end{remark}

Remark \ref{split} gives us 

\begin{remark}\label{splitchi}
	If $\wp $ is a finite 
	place of $K$ that is
	split in the extension $L/K$.
	If the $\wp $-adic Schur index of $\chi $
	is 1, 
	then $\wp $ is not 
ramified in $\Delta (\chi )$.
\end{remark}

\subsection{Modular reduction} \label{decomposition}

Let $\chi \in \Irr^o(G)$.

\begin{theorem}\label{duallatconst}
Let $\wp $ be a finite place of $K$ that is inert in $L/K$.
If $\chi \pmod{\wp }$ is unitary stable 
then $\wp $ is not ramified in the discriminant algebra $\Delta(\chi )$.
\end{theorem}

\begin{proof}
	The proof of this theorem is almost the same as the one 
	of \cite[Theorem 6.4]{OrthogonalStability}, showing 
	the result for real characters $\chi $. 
	As $\wp $ is inert, 
	the completion $L_{\wp }$ of  $L$ at the finite place $\wp $ of $K$
	is  the unique unramified extension of degree $2$ of $K_{\wp}$.
	Denote by $P$ the maximal ideal of 
 $ \Z_{L_{\wp}}$ and by $F$ the residue field $\Z_{L_{\wp}}/P $.
	After possibly passing to a suitable unramified extension of $K_{\wp}$
	we may assume that 
	there is an   $L_{\wp } G$-module $(V,H)$ affording 
	the character $\chi $. 
	Let $\Lambda \subseteq \Lambda^* $ be a maximal 
	integral $G$-invariant lattice in $V$. 
	Then 
	$$P \Lambda ^* \subseteq \Lambda \subseteq \Lambda ^* $$
	and $(\Lambda ^*/\Lambda ,\tilde{H})$ is a Hermitian $FG$-module 
	by Proposition \ref{duallat} (e). 
	Now  $\Lambda $ is maximal integral, so this $FG$-module is the 
	orthogonal 
	direct sum of simple unitary modules. 
	By assumption all these modules have even dimension, 
	so $\dim _{F} (\Lambda ^*/\Lambda )$ is even 
	and Proposition \ref{duallat} (d) shows that 
	the $\wp $-adic valuation of the discriminant of $H$ is even.
\end{proof}

\begin{corollary}\label{ramdiv}
	All prime ideals of $K$ that ramify in the 
	discriminant algebra $\Delta(\chi )$  do
	divide  $2|G|$.
\end{corollary}

For primes that are ramified in $L/K$ the module $\Lambda ^*/\Lambda $ 
carries a symplectic $G$-invariant form. So 
Proposition \ref{duallat} (c)   implies the following result.

\begin{corollary}\label{ramifiedchi}
	Let $\wp $ be a  non-dyadic prime of $K$ that
	is ramified in $L/K$.
	If $\chi \pmod{\wp }$ is 
	unitary stable, then  $\disc(\chi) \pmod{\wp} $ is 
	the orthogonal discriminant of $\chi \pmod{\wp }$.
	This implies that $\wp $ is ramified in the discriminant 
	algebra $\Delta (\chi )$ if and only if the orthogonal discriminant 
	of $\chi \pmod{\wp }$ is not a square in the residue field.
\end{corollary}

For blocks of defect $1$, also the converse of Theorem \ref{duallatconst}
is true. For general cyclic defect 
an analogous statement as \cite[Theorem 6.12]{OrthogonalStability} 
is also true for unitary stable characters. 

\begin{proposition}\label{def1}
	Assume that $\chi $ is in a $p$-block of defect 1
	and let $\wp $ be a finite place of $K$ 
that contains the rational prime $p$
and is inert in $L/K$.
Then $\chi \pmod{\wp }$ is unitary stable 
if and only if $\wp $ does not ramify in the discriminant
algebra of $\chi $.
\end{proposition}

The proof is almost the same as for orthogonally stable characters in
\cite[Theorem 6.10]{OrthogonalStability}. 
In particular this shows that if $\wp $ is inert and $\chi $ 
in a block of defect 1, then
the reduction $\chi \pmod{\wp}$ has at most two unitary 
constituents.
Note that this does not imply that $\chi \pmod{\wp }$ has only 
two constituents, as shown in the following example.

\begin{example}\label{ExON}
As an example let $\chi $ be one of the two complex conjugate 
indicator $o$ irreducible 
ordinary faithful character of degree $116622$ of the group $3.ON$, 
the Schur cover of the sporadic simple O'Nan group of order
$2^9\cdot 3^5\cdot 5\cdot7^3\cdot 11\cdot19\cdot 31$. 
Then $\Q(\chi ) = \Q[\sqrt{-3}]$ and $7,19$ and $31$ are norms 
in $\Q[\sqrt{-3}]/\Q $.
The 11-modular reduction of $\chi $ is irreducible. 
Modulo $5$ the character has defect $1$ and 
$\chi \pmod{5} = 5643ab+52668ab $. 
Whereas $5643a$ and $5643b$ are unitary, the other two 5-modular 
constituents do not admit a non-zero invariant form. 
So Proposition \ref{def1} implies that $5$ divides the discriminant of $\chi $ 
and hence $\disc(\chi ) \in \{ -5 ,-10 \}$. 
Also the 3-modular reduction 
$\chi \pmod{3} =  6138+104346 $ is unitary stable. 
Both 3-modular constituents have indicator $+$, 
$\disc(6138) = O-$ and $\disc(104346) = O+$. 
So the unitary discriminant of $\chi $ is not a square modulo $3$, 
and hence $\disc(\chi ) = -10$. 
\end{example}

\subsection{The unitary discriminants for $O_{10}^+(2)$.}

The group $O_{10}^+(2)$ is 
one example where modular reduction is enough to obtain all 
unitary discriminants.

\renewcommand{\arraystretch}{1.2}
\begin{theorem}
	The unitary discriminants of the characters $\chi $ in 
	$\Irr^o(O_{10}^+(2))$ 
	are as follows: 
$$
\begin{array}{|c|c|c|c|c|}
        \hline
 \chi &  \chi(1) & \Q(\chi) & \disc(\chi) & \Delta (\chi)  \\
        \hline
	33,34 & 110670 & \Q[\sqrt{-15} ] & -1 & [(-1,-3)_{\Q}] \\
        \hline
	51,52 & 332010 & \Q[\sqrt{-15} ] & -2 & [(-2,-5)_{\Q}] \\
        \hline
	68,69 & 442680 & \Q[\sqrt{-15} ] & 1 & [\Q] \\
        \hline
	79,80 & 711450 & \Q[\sqrt{-7} ] & -3 & [(-1,-3)_{\Q}] \\
        \hline
	81,82 & 711450 & \Q[\sqrt{-7} ] & -3 & [(-1,-3)_{\Q}] \\
        \hline
\end{array}
$$
\end{theorem}

\begin{proof}
	The first three pairs of characters $\chi $ 
	with character field $L:=\Q[\sqrt{-15}]$ have a unitary stable reduction modulo all primes that are 
	inert in $L/\Q$. 
	So by Theorem \ref{duallatconst} no inert prime ramifies in 
	the discriminant algebra of $\chi $. 
	Remark \ref{splitchi} tells us that no split prime 
	ramifies in $\Delta (\chi )$ and Remark \ref{infinitechi} 
	shows that the infinite place of $\Q $ ramifies in $\Delta(\chi )$ 
	for the first two but not for the third of the three pairs of 
	characters. 
	For all three characters the 5-modular reduction is  
	irreducible. From the database of orthogonal discriminants in 
	\cite{OD} we get that this 5-modular reduction has orthogonal
	discriminant $O+$ for $\chi _{33/34}$ and $\chi_{68/69}$ and 
	$O-$ for $\chi_{51/52}$. 
	So Corollary \ref{ramifiedchi} shows that $5$ is ramified
	in $\Delta (\chi )$ only for $\chi=\chi_{51/52}$. 
	As the number of ramified places in $\Delta (\chi )$ is even, 
	this yields the behaviour of the prime $3$.
	So 
	$\Delta (\chi _{33})$ is ramified at $\infty $ and $3$, 
	$\Delta (\chi _{51})$ is ramified at $\infty $ and $5$, 
	and 
	$\Delta (\chi _{68})= [\Q ]$.
	\\
	The four characters of degree 711450 are  
	irreducible modulo $5,7,17,$ and $31$. 
	Modulo the ramified prime, 7, their 
	orthogonal discriminant is $O+$ according to \cite{OD}.
	As 2 is a norm in $\Q[\sqrt{-7}]/\Q $ the
	possible discriminants are $-1$ and $-3$. 
	But $-1$ is not a square modulo $7$ and hence 
	$\disc(\chi) = -3$ for these four characters $\chi $.
\end{proof}

\section{Further character theoretic constructions} 

\subsection{Restriction and induction}\label{restind}

Let $G$ be a finite group and $\chi \in \Irr^o(G)$.
Denote by $L = \Q(\chi )$ the character field 
of $\chi $.
If one can construct an $L$-representation $\rho $ affording the 
character $\chi $, then finding a non-degenerate 
$\rho(G)$-invariant Hermitian form boils down to solving a
system of linear equations. A more sophisticated method,
given in \cite{Bernd} to obtain the commuting algebra,
can also be applied for finding invariant Hermitian forms. 

Often the degree of $\chi $ is too big so that we cannot 
construct such a representation $\rho $, however 
sometimes the restriction of $\chi $ to a subgroup remains
unitary stable and all its summands that occur 
in $\chi _{|U}$ with odd multiplicity have manageable degree.
Then we can use Proposition \ref{rechus} to compute $\Delta (\chi )$:

\begin{remark}\label{restriction}
	If there is a subgroup $U \leq G$ such 
	that the restriction $\chi _{|U} $ of $\chi $ 
	to $U$ is unitary stable, then 
	$\Delta (\chi ) = \Delta_L (\chi _{|U})$. 
\end{remark}

Note that  the  parity of the 
 degrees of the unitary constituents  of the restriction 
 matters, but their character fields can be dealt with Corollary \ref{lemnorm}. 
For induction it is the opposite: 
If $\psi $ is a character of some subgroup $U$ and 
$F$ a $U$-invariant form in a representation affording the
character $\psi$, then the orthogonal sum of 
$[G:U]$ copies of $F$ is a $G$-invariant form in 
the induced representation. 

\begin{remark}\label{induz}
	Let $U\leq G$ and let $\psi $ be some character of 
	$U$.
	If the induced character $\psi ^G$ is unitary stable 
	and $[\Q(\psi ) : \Q(\psi ^G)] $ is odd 
	then 
	$$\Delta(\psi ^G) = \left\{ 
	\begin{array}{ll} N_{\Q(\psi ) / \Q(\psi ^G)} (\Delta (\psi )) & \mbox{ if } [G:U] \mbox{ odd } \\
		1 & \mbox{ if } [G:U] \mbox{ even } 
	\end{array} \right.$$
\end{remark}

If the degree $[\Q(\psi ) : \Q(\psi ^G)]$ is even, then we might obtain some 
local information about the discriminant of the induced character by passing to 
completions where this degree is odd. However, in our computations this never gave us
any restriction that we could not obtain by other arguments.

\subsection{Tensor products and symmetrizations}

Tensor products and symmetrizations often give rise to 
equations relating discriminants of irreducible characters:

\begin{remark} \label{tensor} 
	Let $\chi , \psi $ be two  irreducible 
	ordinary characters of $G$ such that $\chi (1)$ is even. 
	If the tensor product $\chi \cdot \psi $ is unitary 
	stable, then $$\Delta_L(\chi \cdot \psi ) = \Delta_L(\chi )^{\psi (1)}$$
	over any field $L$ containing $\Q(\chi )$ and $\Q(\psi )$.
\end{remark}

Symmetrizations have been dealt with in \cite{symmetrisations}
and illustrated with the example of the group $6.Suz$, 
where the symmetrizations of the complex character of degree 12 allow
to conclude that all faithful 
$\chi \in \Irr^o(6.Suz)$ have unitary discriminant $(-1)^{\chi(1)/2}$ 
(see \cite[Section 6]{symmetrisations}).
Symmetrizations are quite helpful for groups $G$ for which one of its 
covering groups has a faithful character $\chi $ of fairly small degree. 
As for induced characters 
it is not important that $\chi (1) $ is even 
to use unitary stable symmetrizations of $\chi $ for getting equations 
for unitary discriminants over the character field of $\chi $.

As an example consider the group $U_4(2)$. 
It has three pairs of complex conjugate indicator $o$ 
irreducible characters of even degree, 
$\chi_{5,6}$ of degree 10, 
$\chi_{12,13}$ of degree 30, 
$\chi_{14,15}$ of degree 40. 
Symmetrizing the indicator $o$ characters of degree 5, we obtain 
$\chi_{5,6}$ as the exterior square, i.e. partition $[1,1]$, 
and $\chi _{14,15}$ as the symmetrization of degree 3 corresponding
to the partition $[2,1]$. 
The sum $\chi_{13} + \chi_{14}$ is the symmetrization of degree 6 
for the partition $[4,1,1]$ of a suitable character of degree 4 of $2.U_4(2)$. 
Combining the results we obtain 
 $\disc(\chi ) = (-1)^{\chi(1)/2}$ for all $\chi \in \Irr^o(U_4(2))$.

\section{Perfect groups with even center} \label{2G}

Let $G$ be a perfect group with a center of even order
and $\chi \in \Irr^o(G)$ be a faithful character of $G$. 
Put $L:=\Q(\chi )$ and $K:=L^+$ the maximal real subfield 
of the character field $L$.

\begin{proposition} \label{Center2}
Let $\wp $ be a non-dyadic prime of $K$.
\begin{itemize} 
\item[(a)]
	If $\wp $ is unramified in $L/K$ then 
		$\wp $ is unramified in $\Delta(\chi )$.
	\item[(b)] 
		If $\wp $ is ramified in $L/K$ and 
		$d$ is the sum of the dimensions of the 
		orthogonal constituents of $\chi \pmod{\wp }$ 
		then 
		$\wp $ is ramified in $\Delta(\chi )$
		if and only if $(-1) ^{d/2} $  
		is not a square modulo $\wp $.
	\end{itemize}
\end{proposition}

	\begin{proof}
(a) For any non-dyadic prime $\wp $ of $K$ the central element $z\in G$
of order $2$ acts as $-id $ on all $\wp $-modular constituents of $\chi $.
	As $G$ is perfect, $\det(-id) = 1$ and thus
 all these constituents have even degree. Hence 
		$\chi \pmod{\wp }$ is unitary stable. 
		Now (a) follows from Remark \ref{splitchi} and 
		Theorem \ref{duallatconst}. 
		\\
		(b)  
	As in (a) all $\wp $-modular constituents of $\chi $ 
		have even degree.
		Let $P$ be the prime of $L$ such that $P^2=\wp \Z_L$. 
		After passing to an unramified extension of the completion $L_P$ 
		we may assume that the Schur index of $\chi $ over $L_P$ is 
		trivial. So there is a Hermitian $L_PG$-module $(V,H)$ affording 
		the character $\chi $. 
		By Proposition \ref{duallat} there is 
		a $\Z_{L_P} G$-lattice $\Lambda $ in $V$ such that 
		$P\Lambda ^* \subset \Lambda \subset \Lambda ^* $.
		Moreover 
		$(\Lambda / P\Lambda^*,\overline{H}) $  is 
		a non-degenerate symmetric bilinear space 
		 over the residue field
		$\Z_L/P \cong \Z_K/\wp $ and 
		$(\Lambda^*/\Lambda ,\tilde{H}) $ is a symplectic 
		$\Z_L/P\Z_L $-module. 
		In particular all orthogonal $\wp $-modular
		constituents occur with
		even multiplicity in $\Lambda^*/\Lambda $, so 
		$$d\equiv \dim(\Lambda/P\Lambda ^*) \pmod{4} .$$
		Clearly $z$ 
		acts as $-id $ on  $\Lambda/P\Lambda^* $.
		Since $G$ is perfect, the Spinor norm of 
	$-id $ is trivial.
	But the Spinor norm of $-id $ 
		is the determinant of  $\overline{H}$  
		(see for instance \cite[Section 3.1.2]{habil})
		and so also this determinant is trivial and 
		 $\disc(\overline{H}) = (-1) ^{d/2} \in \{ 1, -1 \}$.
		Now Corollary \ref{ramifiedchi} tells us that 
		$\wp $ is ramified in $\Delta(\chi )$, if and only if 
		$(-1) ^{d/2}$ is not a square modulo $\wp $.
	\end{proof}

\begin{corollary} \label{Center4}
Let $G$ be
a perfect group whose center has an order divisible by $4$ and let
$\chi \in \Irr^o(G) $ be faithful. 
	Then only dyadic primes of $K$ might ramify in $\Delta (\chi )$.
\end{corollary}

\begin{proof}
Let $z$ be a central element of order $4$. 
Then $z$ acts as a primitive fourth root of unity on 
all irreducible representations of $G$. Again, the determinant
of such a representation is $1$, so all faithful irreducible
characters have a degree that is a multiple of $4$. 
	As this also holds for $\chi  \in \Irr^o(G)$, Remark \ref{infinitechi} 
	tells us that no 
infinite place of $K$ ramifies in $\Delta(\chi )$. 
Proposition \ref{Center2} shows that the only primes of $K$
that might possibly ramify in $\Delta (\chi )$ are 
the dyadic primes of $K$.
\end{proof}

	%\begin{proposition} \label{Center6}
%For a perfect group $G$ whose center has an order divisible by $6$
        %all faithful $\chi \in \Irr^o(G) $ have
        %unitary discriminant $\disc(\chi ) = 1$.
%\end{proposition}

\subsection{The $Q_8$-trick} 

\begin{remark}
	Let $G$ be a group with cyclic center of even order and 
	let $z\in G$ denote the central element of order 2 in $G$.
	Assume that $G$ contains a subgroup $U\cong Q_8$ with 
	$z\in U$. 
	Let $\chi $ be the unique faithful irreducible character of $U$.
	Then the restriction of any faithful character $\Chi $ of $G$ to $U$ is
	a multiple of of $\chi $ and hence unitary stable. 
	Then Remark \ref{simplesub} and Theorem \ref{ind+-} (b) allow 
	to read off 
	the discriminant algebra of $\Chi $ just from its
	degree and its character field: 
	Put $L:=\Q(\Chi) $ and let $K$ denote its maximal real subfield. 
	Then the discriminant algebra of $\Chi $ is
	$$\Delta(\Chi ) = [(-1,-1)_K]^{\Chi(1)/2}  \in \Br(K).$$
\end{remark}

The assumption that $Q_8$ is contained in such a group occurs 
quite frequently. For instance all covering groups $2.A_n$, for 
$n\geq 4$, contain such a group $Q_8$ 
(as $2.A_4$ does) and hence the $Q_8$-trick 
gives all the unitary discriminants of their faithful unitary characters. 

To illustrate the $Q_8$-trick, let $G=2.S_6(3)$. 
Then all faithful $\chi \in \Irr^o(2.S_6(3))$ 
 have character field 
$L=\Q(\sqrt{-3})$. As $(-1,-1)_{\Q} = (-3,-2) _{\Q} $ the remark 
shows that 
$\disc(\chi ) = (-2)^{\chi(1)/2} $.

\section{Subalgebras from suitable automorphisms} \label{antiadj}

\subsection{The fixed algebra under an automorphism}

Let $G$ be a finite group and $\chi \in \Irr^o(G)$.
Put $2m:= \chi (1)$, $L=\Q(\chi )$,
and assume that there is an $L$-representation
$\rho : G \to \GL(V) $ 
affording the character $\chi $.
Put $K:=L^+$ to denote the real subfield of $L$ and let 
$H:V\times V \to L$ be a $\rho(G)$-invariant $L/K$-Hermitian form on $V$.
Assume that there is $\alpha \in \Aut (G)$ 
with $\alpha ^2 = 1$ such that $\chi \circ \alpha  = \overline{\chi }$. 
Then $\chi + \overline{\chi }$ extends to an  
irreducible character 
$$\Chi = {\textrm{Ind}}_G^{\mathcal G}(\chi ) \mbox{ of the semidirect product }
{\mathcal G}:= G:\langle \alpha \rangle .$$ 

\begin{theorem}\label{fixalpha} 
Let
$A:=\Fix_{\alpha }(\rho ) = \langle \rho (g) + \rho(\alpha (g)) \mid g\in G \rangle _K $
denote the $\alpha $-fixed algebra. 
\begin{itemize}
\item[(a)] Let ${\mathcal V}$ be a $K{\mathcal G}$-module 
	affording the character $2\Chi $ and put ${\mathcal Q}:= \End_{K{\mathcal G}} ({\mathcal V}) $.  \\
	Then $[A] = [{\mathcal Q}] \in \Br_2(L,K)$. 
\item[(b)] The algebra $A$ is invariant under the adjoint involution of $H$,
	i.e. $\iota_H(A) = A$. 
	\item[(c)] 
Assume that 
 the Frobenius-Schur indicator of $\Chi $ is $+$.
		\\
		Then $A$
	is an orthogonal subalgebra of $(\End_L(V),\iota _H)$ and  
		$\disc(H) = \disc_L([{\mathcal Q}])^m \disc(\iota _{|A} )$.
	\item[(d)] 
 Assume that the Frobenius-Schur indicator of $\Chi $ is $-$. \\ Then
		$\disc(H) = \disc_L({\mathcal Q})^m$.
	\end{itemize}
\end{theorem}

\begin{proof}
	(a) 
Write $L=K[\sqrt{-\delta }]$, then
$\End_L(V) = A \oplus \sqrt{-\delta } A$, in particular $LA = \End_L(V)
	\cong L^{2m\times 2m}$.
	So $\dim_K(A) = (2m)^2 $ and $A$ is central simple. 
	\\
	To identify the class of the central simple $K$-algebra $A$
	in the Brauer group, we relate the centraliser of 
	$A$ in $\End_K(V) $ to the Schur indices of $\Chi $, 
	i.e. to the class $[{\mathcal Q}] \in \Br_2(L,K)$. 
	We closely follow \cite[Section 2]{Benson}: \\
	Let $\overline{\phantom{x}} $ denote the non-trivial 
	Galois automorphism of $L/K$.
	Identify $\End_L(V)$ with $L^{2m\times 2m}$ by 
	choosing an $L$-basis $B$ of $V$ and  
	define the $K$-linear endomorphism $\sigma \in \End_K(V)$ 
	by
	$$\sigma (\sum _{i=1}^{2m} a_i B_i ) = \sum _{i=1}^{2m} \overline{a_i} B_i .$$
	For $X\in L^{2m\times 2m}$   the Galois conjugate 
	matrix 
	$\overline{X}  = (\overline{X_{ij}}) $ is the matrix of 
	the endomorphism $\sigma  X \sigma $.
	As $\chi \circ \alpha = \overline{\chi }$, there is some 
	$X \in \GL_{2m}(L)$ unique up to scalars in $L^{\times }$ such that 
	$$ X \rho (g) X^{-1}  = \overline{\rho (\alpha(g) )} \mbox{ for all } g\in G .$$
	Then $\overline{X} X $ commutes with all $\rho (g)$, so 
	$\overline{X} X = \lambda I_{2m}$ for some $\lambda \in L^{\times }$. 
	It is easy to see that $\lambda = \overline{\lambda } \in K^{\times }$ 
	and is well defined up to norms.
	Now  the proof of \cite[Proposition 4.2]{Benson} shows that 
	$[{\mathcal Q}] = [(L,\lambda )_K]$. 
	\\
	As $A = \{ \rho (x) \mid x\in KG, \alpha (x) = x \}  $ 
	we see that  $Xa = \overline{a} X $ 
	for all $a\in A$.
As endomorphisms of $V$ this equation translates into 
$X \cdot a = \sigma \cdot a \cdot \sigma \cdot X $.
	So the composition $\beta := \sigma \cdot X$ 
	is a $K$-linear map on $V$ that 
	commutes with all matrices in $A$ and satisfies
	$$\beta^2 = (\sigma X)^2 = (\sigma X \sigma ) X = \overline{X} X 
	= \lambda .$$
	Also scalar multiplication $\gamma $ by $\sqrt{-\delta }\in L$ 
	commutes with all matrices in $A$. 
	These two $K$-linear endomorphisms of $V$ generate a $4$-dimensional
	subalgebra of $\End_K(V)$ and satisfy the relations 
	$$\gamma ^2 = -\delta, \beta ^2 = \lambda , \beta \gamma 
	= - \gamma  \beta $$ 
	so $\langle \gamma ,\beta \rangle _K \cong {\mathcal Q}$.
	\\
	(b) 
 To see  that the adjoint involution $\iota _H$ restricts 
	to a $K$-linear involution on $\Fix_{\alpha }(\rho )$ 
	note that $\iota _H(\rho (g)) = \rho (g^{-1})$ for all $g\in G$,
	so $\iota_H$ maps the generator 
	$\rho(g) + \rho(\alpha (g))$ to the generator
	$\rho(g^{-1}) + \rho(\alpha (g^{-1})) $, in particular 
	$\iota_H(A) = A$. 
\\
	(c)
	To apply Proposition \ref{ind+-} (a) it
 remains to show that 
	the restriction of $\iota _H$ to $A$ is an orthogonal involution 
	on $A$, i.e. that the dimension of skew-adjoint elements in 
	$A$ is $m(2m-1)$. This follows from the assumption 
	on the indicator of $\Chi $. 
	As dimensions remain unchanged under field extensions we 
	may replace $K$ by the field $\R $ of real numbers and
	$L$ by $\C $.
	The assumption on the Frobenius-Schur indicator of $\Chi $ 
	shows that the Schur index of $\Chi $ over the real numbers 
	$\R $ is 1. So
the main result of \cite{Benson} shows that 
there is a basis of $V \otimes _L \C $
	 such that the 
	matrices in $\rho (G)$ with respect to 
	this basis satisfy 
	$$\overline{\rho(g)} = \rho (\alpha (g) ) \mbox{ for all } g\in G.$$
	In particular the set 
	$\rho (G)$ is invariant under complex conjugation and 
	$$A\otimes _K \R \cong \R ^{2m \times 2m} \leq \End_\C (V\otimes _L \C) = \C^{2m\times 2m} .$$
	Also the one dimensional 
	$\R $-space of invariant Hermitian forms is invariant
	under complex conjugation and hence all these forms are real and 
	symmetric, inducing an orthogonal involution on $A\otimes _K \R $.
	\\
	(d) Now assume that the indicator of $\Chi $ is $-$. 
	Then by (a) the algebra $A \cong {\mathcal Q}^{m\times m}$ and 
	the restriction of $\iota $ to $A$ is a symplectic involution.
From  Proposition \ref{ind+-} (b) one gets that $\Delta(H) = [{\mathcal Q}]^m$.
\end{proof}

\begin{definition}
	In the notation of Theorem \ref{fixalpha} 
	the algebra $A$ is called the {\em $\alpha $-fixed algebra}
	of $\chi \in \Irr^o(G)$, the square class 
	$$\disc(\iota _{|A} ) =: \disc^{\alpha }(\chi ) 
	=d (K^{\times})^2 \in K^{\times }/(K^{\times })^2 $$ is called the 
	{\em $\alpha $-discriminant} of $\chi $ 
	and the \'etale $K$-algebra 
	$${\mathcal D}^{\alpha }(\chi ) := K[X]/(X^2-d ) $$ 
	is called the 
	{\em $\alpha $-discriminant algebra} of $\chi $.
\end{definition}

\begin{proposition} \label{Teiler} 
	In the notation of Theorem \ref{fixalpha} 
	all prime ideals of $K$ that ramify in 
	the $\alpha $-discriminant algebra of $\chi $ 
	divide $2|G|$.
\end{proposition}

\begin{proof}
Let $\wp $ be a prime ideal of $K$ that does not divide $2|G|$.
Then $\wp $ is not ramified in $L/K$. 
Passing to the completions at $\wp $ put 
$$\Gamma _{\wp} := \langle \rho(g) \mid g\in G \rangle_{\Z _{K_{\wp}}}\leq {\mathcal A}
\mbox{ and } 
\Delta_{\wp } := \langle \rho (g) + \rho(\alpha (g)) \mid g\in G \rangle _{\Z_{K_{\wp}} } \leq A.
	$$ 
	Because $\wp $ does not divide $|G|$ the
	order $\Gamma _{\wp }$ is a maximal order in ${\mathcal A}$. 
	As $\wp $ does not divide $2 \disc(L/K)$ we have 
	$ \Z_L \Delta _{\wp } = \Gamma _{\wp }$ and
	hence $\Delta _{\wp} $ is a maximal order in $A$. 
	Both orders are invariant under the adjoint involution $\iota _H$.
	In particular  for any $\Delta _{\wp }$-lattice 
	$\Lambda \leq V_{K,\wp}$ also its $H$-dual lattice $\Lambda ^*$
	is $\Delta _{\wp}$-invariant. As $\Delta _{\wp }$ is 
	a maximal order, there is $n\in \Z$ such that 
	$\Lambda ^* = \wp^n \Lambda $ and hence $\wp $ is not 
	ramified in the discriminant algebra of $H$. 
\end{proof}

\section{Orthogonal and unitary condensation}\label{condensation}

\subsection{Orthogonal condensation}\label{orthcond}

Condensation methods play an important role 
in the construction of (modular) character tables and 
decomposition matrices for large finite groups (see for instance \cite{Ryba}
for a brief description of the general idea and one of the many 
papers citing this article for further applications). 
Fixed point condensation is an important tool to 
compute orthogonal discriminants \cite[Section 3.3.2]{survey} 
for large degree characters 
 due to the availability of very sophisticated programs
to compute in large permutation representations (in particular in \cite{MAGMA}). 
Let $U\leq G$ and $W\cong \Z ^{[G:U]}$ and 
$\rho_W: G \to \GL(W) $ be such a permutation representation
of the finite group $G$ on the cosets of the subgroup $U$. 
For $S\leq G$ put 
$$e_S:=\frac{1}{|S|} \sum_{h\in S} h \in \Q G $$ 
to denote the projection on the $S$-fixed points.
Then $e_S$ is an idempotent in $\Z [\frac{1}{|S|}] G$ that 
is invariant under the natural involution $\iota $ of the group algebra
with $\iota (g) = g^{-1}$ for all $g\in G$. 

\begin{remark}
	Let $K$ be a field such that $|S| \neq 0 $ in $K$. 
	Put  $W':= KW \rho_W(e_S)$ and 
	$A:= e_S KG e_S $.
	Then the natural involution $\iota$ on $KG$ 
	endows $A$ with an orthogonal involution 
	that we again denote by $\iota $.
	Moreover $W'$ is an $A$-module 
	whose composition factors $C$ are in bijection to 
	the $KG$-composition factors of $KW$ that contain 
	a non-trivial $S$-fixed space. 
\end{remark}

In practice we never can be sure to have found enough generators for
the condensed algebra $A$. So we need to work with a subalgebra $A'$
for which we know the involution $\iota $.
Character theory predicts the dimensions of 
the composition factors of $W'$, so by computing an $A'$-composition
series of $W'$ we can be sure that the 
$A$-composition factors are isomorphic to the $A'$-composition factors 
of $W'$. For such a composition factor $C$, we can compute the 
$A'$-action, and thus the $A$-action on $C$. 

Now let $V$ be an orthogonal composition factor of $KW$ and assume that 
$Q$ is a non-degenerate $G$-invariant quadratic form on $V$. 
As $e_S =\iota(e_S) $ we have an orthogonal decomposition
$(V,Q)=(V_1,Q_1)\perp (V_2,Q_2)$
where $$V_1 = V\rho_W(e_S) \mbox{ and } V_2= V \rho_W(1-e_S) $$
and $\disc(Q) = \disc(Q_1) \disc(Q_2)$.

Let $$N_S := N_G(S) $$
be the normaliser of $S$ in $G$.
Then $V_2$ is a $KN_S$-module and $Q_2$ is an $N_S$-invariant 
quadratic form on $V_2$. 
If this module $V_2$ is orthogonally stable as $N_S$-module,
then we can predict the orthogonal discriminant of $(V_2,Q_2)$
and obtain the discriminant of $(V,Q)$ by computing the 
one of the condensed module $(V_1,Q_1)$. 
If the restriction of $V$ to $N_S$ is not orthogonally stable
then also $(V_1,Q_1)$ is not an orthogonally stable $N_S$-module
and there is no group that we can use to construct the
invariant quadratic form $Q_1$. 
Thanks to \cite[Proposition 2.2]{orthdet} 
(see Section \ref{antiadj}) we can use the involution 
$\iota $ on $A'$ to compute the discriminant of the
polarisation of $Q_1$:

\begin{remark} (see Remark \ref{discInvo})
	Assume that there is 
	$a\in \rho_{V_1}(A)$ such that 
	\begin{itemize}
		\item[(a)] $\iota(a) = -a $ 
		\item[(b)] $a \in \GL (V_1) $
	\end{itemize} 
	then the determinant of $Q_1$ is 
	$\det(a) (K^{\times })^2$.
\end{remark}

\begin{remark}
	If $K$ is a number field, then the computation of the
	composition factors of the condensed module $Ve_S$ 
	is in general unfeasible as meat-axe methods only work 
	well over finite fields. 
	However, given an a priori list of possible determinants 
	of the composition factor $(V_1,Q_1)$ (e.g. all square 
	classes $d(K^{\times})^2$ 
	for which $K[\sqrt{d}]/K$ is only ramified at primes dividing 
	the group order),
 then it is possible to deduce the correct square class by 
	computing the determinant of $(V_1,Q_1)$ modulo enough 
	well chosen prime ideals (that usually do not divide the group 
	order).
\end{remark}

\subsection{Unitary condensation}\label{unitarycond}

Thanks to Theorem \ref{fixalpha} and Proposition 
\ref{Teiler} we can apply orthogonal condensation
to determine the $\alpha $-discriminant of an indicator $o$ character
and thus the discriminants of the invariant Hermitian forms. 

We aim to compute the unitary discriminant of some 
$\chi \in \Irr^o(G)$.
So assume that we are in the situation of Theorem \ref{fixalpha},
in particular there is some automorphism $\alpha $ of order 2
of $G$ such that $\overline{\chi } = \chi \circ \alpha $. 
Let $V$ be an $LG$-module affording the character $\chi $. 

As before we have to  choose two subgroups of $G$: 
\\
(a) A subgroup $U\leq G$ such that $\chi $ is a constituent of 
the permutation character $\chi _W$ of $G$ on the cosets of $U$. 
\\
(b) A subgroup $S\leq G$ with normalizer $N_S:=N_G(S)$ such that 
\begin{itemize}
	\item[(b1)]
$|S|$ is invertible in $L$,
\item[(b2)]
$V(1-e_S)$ is a 
		unitary stable  $LN_S$-module
	\item[(b3)] $\alpha (S) = S $
\end{itemize}

Because of condition (b3) the automorphism $\alpha $ satisfies
$\alpha (e_S) = e_S$ and hence induces an algebra automorphism 
$\alpha $ on the condensed algebra $A=e_S \rho_W (LG) e_S $. 

To compute the $\alpha $-discriminant of the composition factor 
$Ve_S$ of the natural $A$-module that corresponds to $V$ we compute the
determinant on this composition factor of a 
skew-adjoint unit $X$ in the $\alpha $-fixed algebra in $A$. 
As this singles out a square class of $K^{\times }$,
for which there are only finitely many possibilities thanks to 
Proposition \ref{Teiler}, we can compute composition factors 
and also the determinant of $X$ modulo enough primes 
(not dividing the group order) to deduce the
$\alpha $-discriminant of $Ve_S$. 
The unitary discriminant of
$V(1-e_S)$ can be computed in $N_G(S)$.

\subsection{An example: The group $HN$.}

This example is intended to illustrate the power of condensation methods to handle 
large degree representations. We give some details that should allow the interesting 
reader to develop own condensation programs. 

The sporadic simple Harada-Norton group $HN$ has order 
$2^{14}\cdot 3^6\cdot 5^6 \cdot 7 \cdot 11 \cdot 13 $ and 
two pairs of complex conjugate indicator $o$  irreducible characters 
of even degree: 
$\chi_{25} = \overline{\chi_{26}}$ of degree $656250$ and 
character field $\Q(\chi_{25}) = \Q [\sqrt{-19}]$ and 
$\chi_{35} = \overline{\chi_{36}}$ of degree $1361920$ and 
character field $\Q(\chi_{35}) = \Q [\sqrt{-10 }]$. 
All Schur indices of $HN:2$ are 1. 
\\
$\chi_{25}$ is irreducible modulo $5,7,19$, 
and has two odd degree constituents modulo $11$. 
As 19 is ramified in the character field, its 19-modular reduction is 
orthogonal, and from \cite{OD}  we obtain 
$\disc(\chi_{25} \pmod {19} ) = O+$. 
As $5, 7,$ and $11$ are norms the possible unitary discriminants 
of $\chi _{25} $ are $-2$ and $-3$. 
\\
Similarly $\chi_{35}$ is irreducible modulo $2, 7, 19$ 
and has two odd degree constituents modulo $11$. 
As $11$ is a norm in $\Q[\sqrt{-10}]/\Q $ the possible discriminants 
are $1,3,5,15$.
\\
The missing information is obtained using unitary condensation.
The outer elements in $HN.2$ interchange $\chi_{25} $ and $\chi_{26}$ 
and also $\chi_{35}$ with $\chi_{36}$.
The normaliser in $HN.2$ of the Sylow 5-subgroup does not contain an
outer element of order $2$. So we take a subgroup $S$ of order $5^5$ 
whose normaliser $N_S:=N_{HN} (S)$ is the maximal subgroup of order 
$2000000$ of $HN$. 
Here we find an element $\alpha \in HN.2\setminus HN$ of order 2 
that normalises $S$. 
The group $S$ has a $210$-dimensional fixed space on the 
module with character $\chi _{25}$ and a $416$-dimensional fixed space 
for $\chi_{35}$. 
We start with a permutation representation on the 
$16500000$ cosets of the maximal subgroup $U_3(8).3$ of $HN$. 
All four characters,  $\chi_{25},\chi_{26},\chi_{35},$ and 
$\chi_{36}$ occur with multiplicity one in the corresponding permutation
character.
The group $S$ has $5280$ orbits, so we need to compute 
the action of a skew adjoint $\alpha $-fixed element on 
the composition factors of dimension $210$ respectively 
$416$ of a $5280$-dimensional module.  
We reduce modulo primes where $-19$ (resp. $-10$) are squares. 
It turns out that  the $S$-fixed part 
of the $\alpha $-determinant of both, $\chi _{25}$ and $\chi_{35}$, 
is $33$. 
On the orthogonal complement, the normaliser $N_S$ acts unitary stably 
and we compute the discriminant here as $1$. 
As $11$ is a norm in both imaginary quadratic extensions we get 

\begin{theorem}
The unitary discriminants of the
characters in $\Irr^o(HN)$ are
$\disc(\chi_{25}) = -3$ and $\disc(\chi_{35}) = 3$.
\end{theorem}

\subsection{An example: The group SU$(3,7)$.}

\renewcommand{\arraystretch}{0.8}
I decided to include this example, as the methods generalize to compute 
unitary discriminants of SU$_3(q)$ for arbitrary odd prime powers $q$. 
Let 
$$G:=\mbox{SU}(3,7):=\{ X \in \SL_3(\F _{49})  \mid 
X \left( \begin{array}{@{}c@{}c@{}c@{}} 0 & 0 & 1 \\ 0 & 1 & 0 \\ 1 & 0 & 0 
\end{array} \right) \overline{X}^{tr} = \left( \begin{array}{@{}c@{}c@{}c@{}} 0 & 0 & 1 \\ 0 & 1 & 0 \\ 1 & 0 & 0 
\end{array} \right) \} $$
of order $2^7\cdot 3 \cdot 7^3 \cdot 43 $. 
The characters $\Chi \in \Irr^o(G)$ 
together with their character fields and 
unitary discriminants are given in the following table: 

\renewcommand{\arraystretch}{1.2}

$$
\begin{array}{|c|c|c|c|}
	\hline
	 \Chi & \Chi(1) & \Q(\Chi ) & \disc (\Chi ) \\ 
	 \hline 
	 13,14 & 258 & \Q[\sqrt{-1}] & -7 \\ 
	 \hline 
	 15,16 & 258 & \Q[\sqrt{-2}] & -7 \\ 
	 \hline 
	 27,28 & 344 & \Q[\sqrt{-1}] & 21 \\ 
	 \hline 
	 29-32 & 344 & \Q[\zeta_8] & 7 \\ 
	 \hline 
	 33-36 & 344 & \Q[\zeta_8] & 1 \\ 
	 \hline 
	 37-44 & 344 & \Q[y''_{48}] & 1 \\ 
	 \hline 
	 45-58 & 384 & \Q[u_{43},\sqrt{-43}] & 1 \\ 
	 \hline
\end{array}
$$

A Sylow 7-subgroup $S$ of $G$ consists of the upper triangular 
matrices with 1 on the diagonal; the normaliser $B$ of $S$ is 
the group of all upper triangular matrices 
$B\cong 7^{1+2} : C_{48} $. 
The group $B$ has 48 linear characters, the ones that restrict trivially 
to $S$, one character, $\psi $, of degree $48$, restricting to the sum of the 
other 48 degree one characters of $S$, and $8$ characters 
$\chi_1,\ldots , \chi_8$ of degree 42. 
$\chi _1$ has indicator $-$, $\chi _2$ indicator + and rational character field,
$\chi_3=\overline{\chi_4}$ have character field $\Q(\sqrt{-1})$ and 
$\chi_5,\ldots , \chi_8 $ have character field $\Q(\zeta_8)$. 
By Proposition \ref{ind+-} (b)  the unitary discriminant of $\chi _1$ over the 
relevant character fields is $-7$. 
The characters  $\psi $ and 
$\chi _2$ are orthogonal, so \cite[Corollary 4.4]{Albanian} 
implies that $\disc(\chi _2) = -7$ and $\disc(\psi ) = 1$. 
We compute the remaining unitary discriminants by explicitly constructing 
the representation and an invariant Hermitian form 
over the character field.
We obtain 
$\disc(\chi_3)=\disc(\chi_4) = -7$ and 
$\disc(\chi_5) = \ldots = \disc(\chi_8) = -1$. 

The characters $\Chi _{13} ,\ldots , \Chi_{16} $ have a unitary stable
restriction to $B$ of discriminant $-7$. 
The restriction to $B$ of  $\Chi_{45} , \ldots, \Chi_{58} $ is 
the sum $\psi + \chi_1+\ldots + \chi _8$ and hence unitary stable 
of discriminant $1$. 

The other characters do not restrict unitary stably to $B$; the 
non-unitary stable part consists of the respective $S$-fixed points.
Write $B=S:\langle t \rangle $ where $t$ is some diagonal matrix 
of order $48$. 
One computes with GAP that the unitary stable part $V(1-e_S)$ restricts 
to $B$ with character as in the following table:

$$\begin{array}{|c|c|c|c|c|} 
	\hline
	\Chi_V & \Chi_{V(1-e_S)} &\disc(V(1-e_S)) &\Rho_{Ve_S}(t) & \disc(Ve_S) \\ 
	\hline
	27 & \chi_1+\chi_{3,4}+\chi_{5-8}+\psi  & -7 & (\zeta_{12}^{-1},\zeta_{12}^7) & -3 \\
	\hline
	28 & \chi_1+\chi_{3,4}+\chi_{5-8}+\psi  & -7 & (-\zeta_{12}^{-1},-\zeta_{12}^7) & -3 \\
	\hline
	29 & \chi_1+\chi_2 +\chi_4 + \chi_{5-8} + \psi & -7 & (-\zeta_{24}^{-5},-\zeta_{24}^{11} )  & -3 \\
	\hline
	30 & \chi_1+\chi_2 +\chi_3 + \chi_{5-8} + \psi & -7 & (\zeta_{24},\zeta_{24}^{-7})   & -3 \\
	\hline
	31 & \chi_1+\chi_2 +\chi_4 + \chi_{5-8} + \psi& -7  & (\zeta_{24}^{-5},\zeta_{24}^{11})  & -3 \\
	\hline
	32 & \chi_1+\chi_2 +\chi_3 + \chi_{5-8} + \psi & -7 & (-\zeta_{24},-\zeta_{24}^{-7})   & -3 \\
	\hline
	33 & \chi_{1-4} + \chi_6+\chi_7+\chi_8 + \psi & -1 & (-\zeta_{16},\zeta_{16}) & -2-\sqrt{2}  \\
	\hline
	34 & \chi_{1-4} + \chi_5+\chi_6+\chi_7 + \psi & -1 & (-\zeta_{16}^7,\zeta_{16}^7) & -2-\sqrt{2}  \\
	\hline
	35 & \chi_{1-4} + \chi_5+\chi_7+\chi_8 + \psi & -1 & (-\zeta_{16}^5,\zeta_{16}^5)  & -2+\sqrt{2} \\
	\hline
	36 & \chi_{1-4} + \chi_5+\chi_6+\chi_8 + \psi & -1 & (-\zeta_{16}^3,\zeta_{16}^3)  & -2+\sqrt{2} \\
	\hline
	37 & \chi_{1-4} + \chi_5+\chi_6+\chi_8 + \psi & -1 & (\zeta_{48}^{-7} ,\zeta_{48}) & -3\\
	\hline
\end{array}  $$
All the characters $\Chi_{37-44} $ restrict to $B$ 
as the sum $\chi_1+\chi_2+\chi_3+\chi_4 + \psi + $ three characters of
$\{\chi_5,\ldots , \chi_8 \}$ and the sum of two faithful characters 
of $B/S$. 
This yields that the discriminant of $V(1-e_S)$ is $-1$ for 
$\Chi_{37-44}$. 

To obtain the discriminant of $Ve_S$ just from the action of $t$, 
we consider the automorphism $\alpha $ of $G=U_3(7)$ that is given by 
applying the Frobenius automorphism to the entries of the 
matrices. 
Then $\alpha (t) = t^7$ and $t^{8} $ is in the fixed space of $\alpha $. 

For $\Chi_i$ with $27\leq i \leq 32$ or $37\leq i \leq 44$ the element 
$t^{16}$ acts as
$$\rho _{Ve_S} (t^{16}) = \diag (\zeta _3,\zeta_3^{-1}) $$
on $Ve_S$ so $\det (\rho _{Ve_S}(t^{16}-t^{-16}) ) = -3$ 
in all these cases. Note that $3$ is a norm in $\Q(\Chi_i)/\Q(\Chi_i)^+$ 
for  $i\in\{29,\ldots ,32, 37,\ldots , 44 \}$.

For the four characters $\Chi_{33-36}$ we compute the action of 
the skew adjoint $\alpha $-fixed endomorphism given by
$t+t^7-t^{-1}-t^{-7} $ as 
$$a:=\rho_{Ve_S}(t+t^7-t^{-1}-t^{-7})=\diag(\zeta + \zeta ^7 -\zeta ^{-1} -\zeta ^{-7} , 
-\zeta - \zeta ^7 +\zeta ^{-1} +\zeta ^{-7} ) $$
where $\zeta $ is some primitive 16th root of unity of determinant 
$\det(a) = 8\pm 4 \sqrt{2}$, which is a norm in $\Q(\Chi_i)/\Q(\Chi_i)^+$ for 
$i\in\{29,\ldots, 32\}$.

\section{An example: The sporadic simple O'Nan group}\label{ONan} 

\begin{theorem}
The character fields, unitary discriminants, and discriminant 
algebras of the characters $\chi \in \Irr^o(3.ON)$ 
are given in the following table. 
\end{theorem}

\renewcommand{\arraystretch}{1.3}

\small

$$
\begin{array}{|c|c|c|c|c|}
	\hline
 \chi &  \chi(1) & \Q(\chi) & \disc(\chi) & \Delta (\chi)  \\
	\hline
3,4 & 13376 & \Q[\sqrt{-31} ] & 1 & [\Q ] \\ 
	\hline
	5,6 & 25916 & \Q[\sqrt{-5} ] & 1  & [\Q] \\
	\hline
	31-34 & 342 & \Q[\sqrt{2},\sqrt{-3}] & -1 & [(-1,-1)_{\Q[\sqrt{2}]}] \\
	\hline
	45-48 & 52668 &  \Q[\sqrt{7},\sqrt{-3}] & 8+3\sqrt{7}  & [(-3,8+3\sqrt{7})_{\Q[\sqrt{7}]} ]    \\
	\hline
	53,54 & 63612 & \Q[\sqrt{-3}] & 55 & [(-3,55)_{\Q}] \\ 
	\hline
	57,58 & 116622 & \Q[\sqrt{-3}] & -10 & [(-3,-10)_{\Q}] \\ 
	\hline
59,60 & 122760 &  \Q[\sqrt{-3}] & 1 & [\Q ] \\ 
	\hline
	61-64 & 169290 &  \Q[\sqrt{2},\sqrt{-3}] & -1 & [(-1,-1)_{\Q[\sqrt{2}]}] \\
	\hline
	65,68 & 169632 & \Q[\sqrt{15},\sqrt{-3}]  & -2\sqrt{15} + 15 & [(-3,-2\sqrt{15} + 15 )_{\Q[\sqrt{15}]}] \\
	66,67 & 169632 & \Q[\sqrt{15},\sqrt{-3}]  & 2\sqrt{15} + 15 & [(-3,2\sqrt{15} + 15 )_{\Q[\sqrt{15}]}] \\
	\hline
	69,70 & 175770 & \Q[\sqrt{-3}] & -11 & [(-3,-11)_{\Q}] \\ 
\hline
	71,72 & 207360 & \Q[c_{19},\sqrt{-3}] &  3c_{19}^2 + 10c_{19} + 9 & [
		(-3,\disc(\chi))_{\Q[c_{19}]}]\\
	75,76 & 207360 & \Q[c_{19},\sqrt{-3}] & -13c_{19}^2 + 3c_{19} + 76  & [(-3,\disc(\chi))_{\Q[c_{19}]}]\\
	73,74 & 207360 & \Q[c_{19},\sqrt{-3}] & 10c_{19}^2 - 13c_{19} - 29 & [(-3,\disc(\chi))_{\Q[c_{19}]}]\\
\hline
	77-80 & 253440 & \Q[\sqrt{93},\sqrt{-3}] & 1 & [\Q[\sqrt{93}]] \\
	\hline
\end{array}
$$

\normalsize

\begin{proof}
Note that we only need to treat one of each pair 
	of complex conjugate characters,
	as $\Delta (\chi ) = \Delta(\overline{\chi })$ 
	is the Clifford invariant of  the orthogonal
	character $\chi + \overline{\chi}$. 
\\
\knubbel
	The characters $\chi_{3}, \chi_{4} = \overline{\chi_3} $ 
	are non-faithful. Their character field, $L=\Q[\sqrt{-31}] $, 
	has class number 3 and
	is only ramified at the rational prime $31$. 
	The $p$-modular reduction of $\chi _3$ is unitary stable
        for all primes but $p=7$. 
	 As $7$ is decomposed in
        $L/\Q $, Theorem \ref{duallatconst} yields that $31$ is the
        only prime that possibly ramifies
	in the discriminant algebra $\Delta (\chi _3)$.
	Modulo the ramified prime, 31, 
	the orthogonal discriminant is a square. 
	So Corollary \ref{ramifiedchi} shows that $31$ is not
	ramified in $\Delta (\chi_3)$ and hence $\Delta(\chi_3) = [\Q ]$.
	\\
	\knubbel Also $\chi_5, \chi_6 = \overline{\chi_5} $ are non-faithful.
	Their character field is $L=\Q[\sqrt{-5}]$, has class number 2, and
	is ramified at 2 and 5. 
	The $p$-modular reduction of these two characters is unitary stable
	for all $p \neq 3,7$. As $3$ and $7$ are decomposed in 
	$L/\Q $, Theorem \ref{duallatconst} yields that the 
	only primes that are possibly ramified 
 in the discriminant algebra are $2$ and $5$. 
	The orthogonal discriminant of the $5$-modular reduction of 
	$\chi _5$ is a square mod 5, 
	Corollary \ref{ramifiedchi} implies that the prime 5 is not ramified in
	$\Delta(\chi_5)$. As 
	the number of ramified primes is even, the prime 
	$2$ cannot be the only ramified prime, so
	$\Delta(\chi _5) = [\Q ]$.
	\\ 
	All other characters in the table are faithful. 
	\\
	\knubbel The characters number $31-34$ restrict  irreducibly
	to the maximal subgroup $C_3\times L_3(7).2$. 
	Their unitary discriminant has been computed by David Schlang
	by restricting further to the normaliser of a 7-Sylow subgroup
	of $L_3(7)$. 
	\\
\knubbel Number 45-48 have character field $L:=\Q[\sqrt{7},\sqrt{-3}]$ 
which is ramified at 3 over its
 maximal real subfield $K=\Q[\sqrt{7}]$. 
 Whereas $L$ has class number 2, all ideals in $K$ are principal ideals.
	$K$ has a totally positive fundamental unit $u = 8+3\sqrt{7}$, which
	is  not a norm in $L/K$. 
The only inert prime in $L/K$ that divides the group order is the one
	dividing 2, here $\chi_{45} \pmod{2}$ is unitary stable, 
	so no inert prime ramifies in the discriminant algebra. 
So it remains to decide whether the two prime ideals $\wp_3, \wp_3'$ 
	of $K$ ramify in $\Delta(\chi_{45})$, i.e. 
	$\disc(\chi_{45}) = u$, or $\disc (\chi_{45}) = 1$ and 
	$\Delta(\chi _{45}) = [K]$.
The induced character Ind$_{3.ON}^{3.ON:2}(\chi_{45} )$ has character
field $K$ and Schur index 2 at the two primes $\wp_3, \wp_3'$ of 
$K$ that divide $3$. 
If $\alpha $ is an outer automorphism of order $2$ in $3.ON:2$, then
the $\alpha $-fixed algebra $A$ is equivalent to 
$(-3,u)_K $ by Theorem \ref{fixalpha}. 
However, as $\chi_{45}(1)$ is a multiple of $4$, Theorem \ref{orthsub}
shows that there is no correction factor because of the Schur indices 
of the induced character. So $\disc(\chi_{45}) $ is represented 
by $\disc^{\alpha }(\chi_{45})$. 
To compute this $\alpha $-discriminant we use condensation with the 
derived subgroup $S$ of the normaliser in $3.ON$ of the 7-Sylow subgroup
of order $2\cdot 7^3$. 
$\chi _{45}$ occurs as a constituent in the permutation representation
on the 41,938,920 cosets of a maximal subgroup of index 57 in 
the subgroup $L_3(7)$ of $3.ON$. 
The group $S$ has $61704$ orbits in this permutation character.
After multiplication with a nontrivial central idempotent we
compute composition factors in a $20568$ condensed module.
The determinant of a skew-adjoint element in the $\alpha $-fixed algebra
of the relevant composition factors is either a square (y) or not (n) 
modulo suitable primes:
$$
\begin{array}{|l|c|c|c|c|c|c|c|c|c|c|}
	\hline
	\chi & \dim &  37 & 43 & 61 & 67 & 73 & 79 & 97 & 103 & 109  \\
	\hline
	 45 & 72 & yn & n & y & y & y & n & n & yn & yy \\
	 \hline
	 57 & 162 & yy & yy & nn & yy & nn & yy & nn & nn & yy \\
	 \hline 
\end{array} 
$$
Note that $7$ is a square mod $37$, $103$, and $109$ but not modulo the
other primes, so one obtains for $\chi _{45}$ two constituents 
of dimension 72 modulo $37$, $103$, and $109$ and one constituent of 
dimension 144 for the other primes.
The results allow to conclude that the $\alpha $-discriminant 
on the $72$-dimensional $S$-fixed space of $\chi_{45}$ is 
$u\cdot 7 \cdot 31$. 
\\
\knubbel Theoretical  arguments showing that
$\disc(\chi_{57}) = -10$ are given in Example \ref{ExON}.
However, during the condensation of $\chi_{45}$, we also computed the
$\alpha $-discriminant of $\chi_{57}$ on the 162 dimensional fixed space of the
$S$
as $\disc^{\alpha}(Ve_S) = -21 $.
This is in accordance with Theorem \ref{orthsub}.
The character field of $\chi_{57}$ is $\Q[\sqrt{-3}]$.
The induced character Ind$_{ON}^{ON:2}(\chi_{57} )$ has character
field $\Q$ and Schur index 2 at the primes 2 and 5.
By Theorem \ref{fixalpha} the $\alpha $-fixed algebra is Brauer
equivalent to $(-3,10)_{\Q }$ and of $\Q[\sqrt{-3}]$-discriminant $10$.
As $3$ and $7$ are both norms in $\Q[\sqrt{-3}]/\Q $ Theorem \ref{orthsub}
yields $\disc (\chi_{57} ) = - 10$.
	\\
	\knubbel The character field of the characters number 53,54,57-60 is
	$\Q[\sqrt{-3}]$ and of class number 1. 
	The prime divisors $7,19,31$ of the group order are split 
	in this extension and hence norms. 
	\\
	\knubbel For the primes 
	5 and 11 the modular reduction of $\chi _{53} $ is not 
	unitary stable. As both primes divide the group order only 
	with multiplicity 1, Proposition \ref{def1} yields
	that 55 divides the unitary discriminant of $\chi _{53}$.
	So $\disc(\chi_{53}) \in \{ 55, 110 \}$. 
	Now $\chi _{53} \pmod{3}$ is orthogonally stable and has 
	square discriminant, so $\disc(\chi_{53}) $ is a square modulo 3,
	allowing to conclude that $\disc(\chi_{53}) = 55$.
	\\
\knubbel The character $\chi_{59}$ reduces irreducibly modulo 5 and 11 
leaving $1$ and $2$ as possible unitary discriminants. Its reduction
modulo 3 is not orthogonally stable, so we cannot obtain any more 
information about this character from the decomposition matrix. 
So we use the unitary condensation method described in Section \ref{unitarycond}. 
We start with the permutation module $P$ of dimension 
$736560$ on the cosets of the subgroup $L_3(7)\leq 3.ON$
and use fixed point condensation with the 7-Sylow subgroup $H$ of $3.ON$. 
The group $H$ has 2196 orbits on the 736560 cosets, 
falling into 732 orbits under the center $\langle z \rangle $ of $3.ON$. 
Computing modulo primes $p \equiv 1 \pmod{3} $ we specify 
the action of $z$ as some primitive third root of unity in $\F_p^{\times} $
and hence reduce the computations to find composition 
factors of a 732 dimensional $\F_p$-module for in total 25 
suitable primes. 
These composition factors have dimension $3,171,192,366$, 
where the latter two correspond to characters number $53/54$ resp. $59/60$. 
For some $\alpha \in N_{3.ON.2}(H) $ of order $2$ that is not 
contained in $3.ON$ we compute the determinants of an $\alpha $-fixed 
skew-adjoint element in the condensed algebra 
 as $3\cdot 5 \cdot 11$ on the composition factor of dimension 192
and $7$ for the one of dimension 366. 
On the orthogonal complement of the fixed space the relevant
discriminants of the restriction to the normaliser of $H$ are
powers of $7$, so one concludes that $\disc(\chi_{59}) = 1$.
\\
\knubbel The characters number 61-64 are tensor products, 
$169632 = 342 \cdot 495$ of one of the characters number 31-34 with 
a character of degree 495. So they have the same unitary discriminant 
as $\chi_{31}$. 
\\
\knubbel Number 65-68 have character field $L:=\Q[\sqrt{15},\sqrt{-3}]$  an 
unramified extension of its maximal real subfield $K=\Q[\sqrt{15}]$. 
Both fields, $L$ and $K$ have class number 2. 
The fundamental unit of $K$ is a norm.
The prime ideals of $K$ that are inert in $L/K$ 
and divide the group order are the divisors of $2$, $5$, $11$, and $19$.
The reduction modulo both prime ideals of $K$ that divide 19 is 
irreducible. 
For 11, one prime ideal yields an  irreducible 
reduction, the other one a non unitary stable reduction. 
By Proposition \ref{def1}  exactly one of these two prime ideals, say 
$\wp _{11}$, ramifies in $\Delta (\chi_{65})$. 
Also the reduction modulo the unique prime ideal $\wp_5$ of $K$ that 
divides $5$ is not unitary stable. 
The ideal $\wp_3\wp_5\wp_{11} = (2\sqrt{15}+15) $ 
has a totally positive generator, but no such combination 
involving the prime ideal dividing 2 has such a totally positive generator. 
So we conclude that the discriminant of $\chi _{65}$ 
is represented by $\pm 2\sqrt{15} + 15$.
\\
\knubbel The character $\chi _{69} = \overline{\chi_{70}}$ has 
character field $\Q[\sqrt{-3}]$, is orthogonally stable modulo 3
with square discriminant. Its reduction modulo
5 is absolutely irreducible and modulo 11 we obtain
$175770 = 42687+133083$. So Proposition \ref{def1}, Theorem \ref{duallatconst}
and Corollary \ref{ramifiedchi} together show that
$\disc(\chi_{69}) = -11 $.
\\
\knubbel The characters number $71-76$ are algebraic
conjugate and so are their unitary discriminants.
It is enough to treat $\chi_{71} = \overline{\chi_{72}}$.
The character field is $L=K[\sqrt{-3}]$ where $K=\Q[c_{19}]$
is the real subfield of $L$.
The extension $L/K$ is only ramified at 3,
$L$ has class number 3, $K$ has class number 1 and all
totally positive units in $K$ are squares.
The prime divisors of the group order that are inert in $L/K$
are $2,5$ and the three prime ideals of $K$ that divide 11.
Whereas $\chi _{71} \pmod{2}$ and $\pmod{5}$ is absolutely
irreducible and hence unitary stable, its $11$-modular reduction
is unitary stable at the prime ideal of $K$ that contains
$c_{19}-6$ and not unitary stable at the other two prime ideals.
So the unitary discriminant is a totally positive generator
of the product of these two prime ideals dividing 11 in $K$.
One checks that this number is indeed a square modulo 3,
as predicted by the orthogonal discriminant of $\chi_{71} \pmod{3}$.
\\
\knubbel The last characters, $\chi_{77},\ldots ,\chi_{80} $, have
character field $L:=\Q[\sqrt{93},\sqrt{-3}]$ again of class number 3.
The fundamental unit of $K=\Q[\sqrt{93}]$ is totally positive
and a norm in $L/K$.
The extension $L/K$ is totally unramified.
The inert primes dividing the group order are 3 and 11.
As both reductions $\chi_{77} $ modulo the two prime
ideals of $K$ that divide 11 remain absolutely irreducible,
the only possible ramified prime in $\Delta (\chi_{77})$ is
the one that divides $3$.
As the number of ramified places is even this prime is also
not ramified and hence the unitary discriminant of
$\chi _{77}$ is 1.
\end{proof}

\end{document}